\documentclass{article}
\title{Invariant forms, associated bundles and Calabi-Yau metrics}
\author{Diego Conti}

\usepackage[english]{babel}
\input xy
\xyoption{all}
\usepackage{amssymb,amsmath,amsthm}

\newcommand{\norm}[1]{\left\Vert#1\right\Vert}
\newcommand{\abs}[1]{\left\vert#1\right\vert}
\newcommand{\R}{\mathbb{R}}

\newcommand{\lie}[1]{\mathfrak{#1}}     
\newcommand{\Lie}{\mathcal{L}}          
\newcommand{\Z}{\mathbb{Z}}

\newcommand{\N}{\mathbb{N}}
\newcommand{\C}{\mathbb{C}}
\newcommand{\CP}{\mathbb{CP}}

\newcommand{\hook}{\lrcorner\,}
\newcommand{\LieG}[1]{\mathrm{#1}}      
\newcommand{\SU}{\mathrm{SU}}

\newcommand{\Gtwo}{\mathrm{G}_2}

\newcommand{\su}{\mathfrak{su}}

\newcommand{\GL}{\mathrm{GL}}

\newcommand{\dfn}[1]{\emph{#1}}
\newcommand{\id}{\mathrm{Id}}   

\newcommand{\tf}[1]{#1}             
\newcommand{\tm}[1]{\overline{#1}}  
\newcommand{\vb}[1]{\underline{#1}} 

\DeclareMathOperator{\Stab}{Stab}
\DeclareMathOperator{\End}{End}

\DeclareMathOperator{\Hom}{Hom}

\DeclareMathOperator{\Diff}{Diff}

\theoremstyle{plain}
\newtheorem{proposition}{Proposition}
\newtheorem{theorem}[proposition]{Theorem}
\newtheorem{lemma}[proposition]{Lemma}
\newtheorem{corollary}[proposition]{Corollary}
\theoremstyle{definition}

\theoremstyle{remark}
\newtheorem*{remark}{Remark}

\theoremstyle{plain}
\newtheorem{algo}[proposition]{Algorithm}
\newcommand{\card}[1]{\abs{#1}}
\DeclareMathOperator{\Span}{Span}

\begin{document}
\maketitle
\begin{abstract}
We develop a method, initially due to Salamon, to compute the space of ``invariant'' forms on an associated bundle $X=P\times_G V$, with a suitable notion of invariance. We determine sufficient conditions for this space to be $d$-closed. We apply our method to the construction of Calabi-Yau metrics on $T\CP^1$ and $T\CP^2$.
\end{abstract}
\vskip5pt\centerline{\small\textbf{MSC classification}: 53C25; 53C30, 57S15, 68W30}

\centerline{\small\textbf{Keywords}: Vector bundle, special geometry, Calabi-Yau, symplectic cone.}

\section*{Introduction}
This paper is motivated by the study of special geometries. A special geometry on a manifold $X$ is a $G$-structure defined by one or more differential forms on $X$, so that
\[G=\{g\in\GL(n,\R)\cong\GL(T_xX)\mid g \text{ preserves the defining forms at $x$}\}\]
for all $x$ in $X$. Having fixed both $G$ and the defining forms, the set of isomorphisms $\GL(n,\R)\cong\GL(T_xX)$ satisfying the above condition defines the $G$-structure as a reduction to $G$ of the bundle of frames. The notion of special geometry also entails  closedness, or more general differential relations among the defining forms. Calabi-Yau manifolds, or more generally manifolds with special holonomy, are examples of special geometries.

In practice, finding explicit examples of special geometries involves solving certain PDE's; to reduce the number of variables, one typically imposes extra symmetry conditions, assuming that a (compact) group of diffeomorphisms $H$ preserves the $G$-structure. Then, about each special orbit, the manifold has the form of an associated bundle
\begin{equation}
\label{eqn:HomogeneousCase}
X=H\times_{H_x}V\;.
\end{equation}
Thus, the defining forms are to be sought amongst invariant forms on an associated bundle. Notice that $H\to  (H/H_x)$ is essentially a $H_x$-structure.

In this paper we develop a method, initially introduced by Salamon in \cite{BryantSalamon,Salamon:RedBook}, to compute the space of invariant forms on an associated bundle
$X=P\times_G V$,
 associated to a generic $G$-structure $P$. Thus, we do not assume that the base manifold $P/G$ is a homogeneous space.

This  method consists in singling out certain local forms on $X$; in particular, one must choose suitable  fibre coordinates $a_i$.
Having fixed a connection, the tangent bundle $TX$ splits into a horizontal and a vertical distribution, and the local one-forms $da_i$ define  forms $b_i$ by projection on the vertical component. Whilst the $a_i$ and $b_i$ are not globally defined in themselves, they can be contracted to obtain globally defined forms, namely
\[\sum_i a_ia_i,\quad \sum_i a_ib_i.\]
Moreover, one can usually find horizontal local forms $\beta_i$ such that
\[\sum_i a_i\beta_i,\quad \sum_i b_i\beta_i,\]
are globally defined forms. In this respect,  the indexed objects $a$, $b$ and $\beta$ resemble the letters of an alphabet, which carry no meaning in themselves, but combine together to produce all the words of a dictionary. In the homogeneous case \eqref{eqn:HomogeneousCase}, it turns out that the elements of the dictionary are invariant under the global action of $H$.

\medskip
The ``dictionary'' method was introduced in \cite{Salamon:RedBook} as a collection of examples, more than a general method, nor was the notion of invariance investigated there. In this paper we give a definition of invariance for forms on $X$ that generalizes the notion of $H$-invariant forms in the homogeneous case \eqref{eqn:HomogeneousCase}. In fact, we call these forms \emph{constant}, since they can be read as constant maps from $P$ to $\Omega(V,\Lambda^*\R^n)^G$; in particular, zero-degree constant forms can be identified with $G$-invariant functions on $V$. We prove that the space of constant forms is finitely generated as a module over the ring of zero-degree constant forms, and that one can generalize the dictionary method  to obtain the full space of invariant forms from a finite number of ``letters''. We actually describe an algorithm to accomplish this in an efficient way, at least in the case where $G$ acts transitively on the sphere in $V$. As an application, we compute the space of invariant forms on $T\CP^2$, and obtain a two-parameter family of local Calabi-Yau metrics with conical K\"ahler form, generalizing the known conical hyperk\"ahler example. This is achieved by considering an associated special geometry in seven dimensions, analogous to Hypo-contact geometry in five-dimensions \cite{DeAndres:HypoContact}, and constructing explicit examples  on the sphere bundle in $T\CP^2$; an extension theorem of \cite{ContiFino} completes the construction.

In order for the space of constant forms $\Omega(V,\Lambda^*\R^n)^G$ to be useful for constructing special geometries, one needs to compute how the operator $d$ acts on its elements. In Lemma~\ref{lemma:Formulazza} we give a formula to this effect. From this formula, we deduce  a sufficient condition for the algebra of constant forms to be $d$-closed, namely that both torsion and curvature be themselves constant. The case~\eqref{eqn:HomogeneousCase} trivially satisfies this condition (although, in this hypothesis, where ``constant'' coincides with ``invariant'', the result is also trivial). Conversely, assuming that $d$ preserves the space of invariant forms, we show that the curvature must be constant, under the additional hypothesis that the infinitesimal action of $G$ on $V$ is effective. We obtain a weaker but similar result for the torsion.

\medskip
The contents of this paper are organized as follows. The first section serves as an illustration of the dictionary method; it contains an elementary example from \cite{Salamon:Explicit} concerning the associated bundle  $T^*S^2$. In particular, a one-parameter family of complete hyperk\"ahler metrics is constructed, giving a different characterization of the Kronheimer metrics on the resolution of $\C^2/\Z_2$ \cite{Kronheimer:ALE}.

In the second section, we show how forms over $X$ can be viewed as sections of a vector bundle associated to $P$, with infinite-dimensional fibre $\Omega(V,\Lambda^*\R^n)$. We introduce the natural action of the gauge group on this space.

Constant forms are introduced in Section~\ref{sec:Closure}; they are characterized as forms that are both gauge-invariant and parallel as sections of the associated bundle introduced in Section~\ref{sec:Gauge}. The problem of whether the algebra of constant forms is $d$-closed is discussed here.

In the fourth section we show that the space of invariant forms is a finitely generated module over the ring of invariant functions on the fibre $V$, and we establish a criterion to determine whether a module of invariant forms coincides with the full space of invariant forms.

The formal definition of the dictionary is introduced in Section~\ref{sec:Dictionary}. This is done under the hypothesis that $M$ is a homogeneous space, with no significant loss of generality, since the results of Section~\ref{sec:Closure} lead one to assume that any two neighbourhoods are affine isomorphic. This section also contains a proof of the fact that one can always find a finite set of letters that generate the full space of invariant forms.

In the sixth section we describe an algorithm, valid in the cohomogeneity one case, to compute the dictionary generated by a certain set of letters, and thus check whether this dictionary is complete.

In the final section we apply the method to $T\CP^2$, finding explicit Calabi-Yau metrics on this space.

\section{Invariant forms on $T^*S^2$}
\label{sec:TS2}
As an introduction to our general construction,  we begin this paper with an elementary example, following \cite{Salamon:Explicit}. Thus, we introduce a list of ``invariant'' forms on $T^*S^2$, and we use this list to obtain a one-parameter family of hyperk\"ahler structures, containing in particular the Eguchi-Hanson metric \cite{EguchiHanson}. We shall not go into detail concerning the relevant notion of invariance, nor give formal descriptions of the objects involved, since this will be the object of the sections to follow.

Let $x_1$, $x_2$ be the coordinates on $S^2$ given by stereographic projection, so that
\[\beta_1=\frac{2}{(x_1)^2+(x_2)^2+1}\,dx_1,\quad \beta_2=\frac{2}{(x_1)^2+(x_2)^2+1}\,dx_2\]
is an orthonormal basis for the standard metric, defined away from the north pole. Then the Levi-Civita covariant derivative gives
\[\nabla \beta_1=\omega\otimes \beta_2,\quad \nabla \beta_2=-\omega\otimes \beta_1,\]
where $\omega$ is the connection form. In this setting, $\omega$ is just a local $1$-form on $S^2$. Since the Levi-Civita connection is torsion-free, it follows that
\[d e^1=\omega\wedge \beta_2,\quad d \beta_2=-\omega\wedge \beta_1\;,\]
Explicitly,
\[\omega=
x_2 \beta_1-x_1 \beta_2\;,\]
and, consistently with the fact that the sphere has curvature $1$,
\[d\omega=-\beta_1\wedge \beta_2.\]

The total space $X$ of $\pi\colon T^*S^2\to S^2$ has a tautological 1-form $\tau$; indeed, a point $x$ of $X$ is a $1$-form on $T_{\pi(x)}S^2$, and one defines
\[\tau_x=\pi^*x\;.\]
We can then introduce fibre coordinates $a_1,a_2$, so that $\tau$ equals \mbox{$a_1\pi^*\beta_1
+a_2\pi^*\beta_2$}.
We shall implicitly pull the $\beta_i$ back to $X$, and write
\begin{equation}
\label{eqn:tau}\tau=\sum_i a_i\beta_i.
\end{equation}
Then
\[
d\tau=\sum_i b_i\wedge \beta_i\quad\text{where}\quad
\left\{
\begin{aligned}
b_1&=da_1-a_2\omega\\
b_2&=da_2+a_1\omega
\end{aligned}
\right.
\]
So at each point we have a basis of $T^*X$, given by $\beta_1$, $\beta_2$, $b_1$, $b_2$. Notice that the annihilator of $\langle b_1,b_2\rangle$ is nothing but the Levi-Civita horizontal distribution, showing that this construction depends on the choice of a connection in an essential way.

Whilst the objects $a_i$, $b_i$ and $\beta_i$ are only defined outside of the north pole, the form $\tau$ is defined globally, and so if of course $d\tau$. In fact, one can obtain a list of globally defined forms on $X$ by contracting the  ``letters'' $a_i$,  $b_i$ and $\beta_i$ as in \eqref{eqn:tau}.
Explicitly, the  ``words'' of our dictionary are
\begin{equation}
\label{eqn:contractions}
\sum_i a_ia_i,\quad \sum_i a_ib_i,\quad \sum_i a_i\beta_i,\quad\sum_i b_i\beta_i\;,
\end{equation}
where we have used componentwise product as a contraction,
\begin{equation}
\label{eqn:contractions2}
a_1b_2-a_2b_1,\quad a_1\beta_2-a_2\beta_1,\quad b_1\beta_2-b_2\beta_1,
\end{equation}
where we have used the determinant as a contraction, and the volume forms on fibre and base
\begin{equation}
\label{eqn:contractions3}
b^1\wedge b^2,\quad \beta_1\wedge \beta_2,
\end{equation}
which can also be  obtained from a contraction using the determinant. The forms \eqref{eqn:contractions} have obvious interpretations: $\sum_i a_ia_i$ is the function $r^2$, given by
\[r^2\colon X\to\R,\quad r^2(x)=\norm{\pi(x)}^2\;,\]
and the others are $rdr$, $\tau$ and $d\tau$. The forms \eqref{eqn:contractions2} arise essentially from the fact that $S^2$ coincides with $\CP^1$, so that one has a natural complex structure both horizontally and vertically.

Having this dictionary of forms at our disposal, we can look for special geometries on $X$. We obtain the following:
\begin{theorem}[\cite{Kronheimer:ALE,Salamon:Explicit}]
The total space of $T^*S^2$ has a one-parameter family of hyperk\"ahler structures, defined by
\begin{align*}
\omega_1&=(k+r^2)^{1/2}\beta_1\wedge \beta_2-(k+r^2)^{-1/2}b_1\wedge b_2,\\
\omega_2&=\beta_1\wedge b_1+\beta_2\wedge b_2,\\
\omega_3&=\beta_1\wedge b_2-\beta_2\wedge b_1,
\end{align*}
where $k$ is a non-negative constant.
\end{theorem}
\begin{proof}
Having fixed coordinates on $S^2$, it is straightforward to check that the $\omega_i$ are closed. By \cite{Hitchin:SelfDuality}, this implies that the $\omega_i$ are parallel and the structure is hyperk\"ahler.
\end{proof}
Notice that for $k=0$ the underlying metric is flat, and degenerate at the zero section; however, it can be made non-degenerate by replacing the zero section with a point. For $k=1$, we obtain the Eguchi-Hanson metric \cite{EguchiHanson}.
\section{Gauge group actions}
\label{sec:Gauge}
Let  $p\colon P\to M$ be a principal bundle with compact fibre $G$, and let $T$ be a representation of $G$ such that
\[P\times_G T \cong TM\;.\]
Although the results of this section do not require this hypothesis, we shall also assume, for simplicity, that $G$ acts effectively on $T$; then $P$ is a reduction of the bundle of frames on $M$, i.e. a $G$-structure.
Now fix another $G$-module $V$, and let $X$ be the total space of the associated
vector bundle
\[\vb{V}=P\times_G V\;.\]
We are ultimately interested in constructing special geometries, by singling out ``canonical'' closed forms on $X$. As  functions on $M$ pull back to functions on $X$, the space of forms $\Omega(X)$ is an algebra over $C^\infty(M)$, clearly not finitely generated. In this section we shall obtain a first reduction of $\Omega(X)$ to a subalgebra that is, on the contrary, finitely generated over $C^\infty(M)$. This subalgebra arises as the space of forms that are invariant under a certain action of the gauge group.

By the definition of  frame bundle, each point of $P$ gives rise to an isomorphism $T\to T_yM$, and so of course $G\subset \GL(T)$ acts on $P$ on the right. However, on each fibre $P_y$ one also has a left action of $G$ given by the inclusion in $\GL(T_yM)$. More precisely, consider the bundle $\vb{G}=P\times_G G$, where $G$ acts on itself by conjugation;
the \dfn{gauge group}  $\mathcal{G}$, as introduced in \cite{AtiyahHitchinSinger},  is the space of its sections, to be denoted by the symbol $\vb{g}$.
Indeed, the bundle map
\begin{align*}
\vb{G}\boxtimes \vb{G}&\to\vb{G}\\
[u,g]\times[u,h]&\to[u,gh]
\end{align*}
induces a group structure on $\mathcal{G}$. Then $\mathcal{G}$ acts on $P$ on the left by
\[\vb{g}\cdot u=ug\;,\quad \text{ where } \vb{g}_{p(u)}=[u,g]\;.\]
In particular, $\mathcal{G}$ acts on every bundle associated to $P$: if $V$ is any $G$-space, i.e. a manifold with a smooth left $G$ action, then the associated bundle $\vb{V}=P\times_G V$ has a left $\mathcal{G}$ action defined by
\begin{equation}
\label{eqn:ActionOnPInducesActionOnX}
\Phi\colon\vb{G}\to\End(\vb{V}),\quad \Phi_{\vb{g}}[u,v]=[\vb{g}\cdot u,v].
\end{equation}
The definition of this induced action only depends on the fact that $\mathcal{G}$ acts on $P$; however, using the definition of $\mathcal{G}$, one can also regard the gauge group action as induced by the bundle map
\begin{equation}
\label{eqn:InducedGaugeAction}
\begin{aligned}
\vb{G}\boxtimes \vb{V}&\to\vb{V}\\
[u,g]\times[u,v]&\to[u,gv]
\end{aligned}
\end{equation}
In particular, the gauge group $\mathcal{G}$ acts on $X$ as a group of diffeomorphisms, and we could let $\mathcal{G}$ act on $\Omega(X)$ accordingly. However, the notion of invariance that we wish to investigate deals with a different action, that reflects the fact that $\mathcal{G}$ acts infinitesimally on $M$, beside acting on each fibre.

\smallskip
Before introducing this action, we need to introduce more language. Let $\pi:X\to M$ be the projection; then the pullback to $X$ of $P$
is \[\pi^*P=P\times V\;,\] with principal $G$ action given by
\[ (p,v)g=(ph,h^{-1}v)\;.\]
By choosing a connection on $P$ we can
make $\pi^*P$ into a $G$-structure on $X$. Indeed, fix a connection form $\omega$ and let $\tf{\theta}$ be the tautological form on $P$, taking values in $T$. Given a vector space $W$, we denote by $\Lambda^*W$ the exterior algebra over $W$. Since $G$ is compact, we can fix  invariant metrics on $T$ and $V$,  identifying each with its dual accordingly; so, we shall make no distinction between, say, $\Lambda^*T$ and $\Lambda^*(T^*)$.
\begin{lemma}
The map
\begin{equation}
\label{eqn:lambda}
\begin{aligned}
T(P\times V)&\to \pi^*P\times (T\oplus V)\\
(u,v;u',v')&\to \left((u,v), \left(\theta_u(u'),-\omega_u(u')v+v'\right)\right)
\end{aligned}
\end{equation}
induces an isomorphism $TX\cong\pi^*P\times_G (T\oplus V)$.
In the same way, one obtains an isomorphism
\begin{equation}
\label{eqn:LambdaX} \Lambda(X)\cong\pi^*P\times_G \Lambda^*(T\oplus V)\;,
\end{equation}
\end{lemma}
\begin{remark}
The identification \eqref{eqn:lambda} depends on the choice of a connection, and so does all of the construction to follow.
\end{remark}
Now let $\mathcal{G}$ act on the left on $\pi^*P=P\times V$ by
\[\Psi\colon\mathcal{G}\to\Diff(\pi^*P),\quad \Psi_{\vb{g}}(u,v)=(ug,v) \text{ if } \vb{g}_{[u,v]}=[u,g]\;.\]
Given a subset $A\subset M$, not necessarily open, let $X_A=\pi^{-1}(A)$. We let
  $\mathcal{G}$ act on $\Gamma(X_A,\pi^*P)$ as follows:
\begin{equation}
\label{eqn:GaugeGroupActsOnPiStarP}
\vb{g}\cdot s = \Psi_{\vb{g}}\left(s\circ\Phi_{\vb{g}^{-1}}\right).
\end{equation}
Unlike its action on $P$, this action of $\mathcal{G}$ is only defined on sections. Therefore, there is no induced $\mathcal{G}$ action on bundles associated to $\pi^*P$, but there is one such action on the space of their sections. Explicitly, we can extend the action $\Psi$ to every associated bundle $\pi^*P\times_G Y$ by
\[\Psi_{\vb{g}}\bigl[s,y\bigr]=\bigl[\Psi_{\vb{g}}s,y\bigr].\]
With this definition, \eqref{eqn:GaugeGroupActsOnPiStarP} also applies to associated bundles. In particular, we have constructed an action of the gauge group on $\Omega(X)$.

\bigskip
We now show that $\Omega(X)$ can be expressed as the space of sections of a vector bundle over $M$ associated to $P$, so that the action of $\mathcal{G}$ on $\Omega(X)$ is the induced gauge action, in the sense of \eqref{eqn:ActionOnPInducesActionOnX}. Let $Y$ be a $G$-space. Then $G$ acts on the left on the space
 \[R=C^\infty(V,Y)\]
 of
$Y$-valued functions on $V$;
the action is given by \[(gr)v=g(r(g^{-1}v))\;.\]
Fix a point $y$ in $M$; then the fibre $X_y$ consists of points $[u,v]$, where $v$ is in $V$ and $u$ is in $P_y$. Consider the map

\begin{equation}
\label{eqn:mu_y}
\begin{aligned}
\mu_y\colon\vb{R}_y&\to\Gamma(X_y,\pi^*P\times_G Y)\\
\mu_y([u,r])_{[u,v]}&=[(u,v),r(v)]
\end{aligned}
\end{equation}
It follows from the definition of the action on $G$ on $R$ that this map does not depend on the choice of $u$. For every $U\subset M$, the $\mu_y$ combine to give a map
\[\mu_U:\Gamma(U,\vb{R})\to\Gamma(X_U,\pi^*P\times_G Y)\;.\]
The gauge group acts on both sides, by (\ref{eqn:InducedGaugeAction}) and \eqref{eqn:GaugeGroupActsOnPiStarP} respectively.
\begin{proposition}
\label{prop:MuIsIsomorphism} The maps $\mu_U$ are $\mathcal{G}$-equivariant isomorphisms.
\end{proposition}
\begin{proof}
It is sufficient to prove the statement for each map $\mu_y$. It is straightforward to show that $\mu_y$ is an isomorphism. Let $\vb{g}$ be an element of the gauge group with $\vb{g}_y=[u,g]$; its action on
$\vb{R}_y$ satisfies
\[\mu_y\bigl(\vb{g}[u,r]\bigr)_{[u,gv]}=\mu_y([ug,r])_{[ug,v]}=\bigl[ug,r(v)\bigr]\;.\]
On the other hand
\[\Psi_{\vb{g}}\mu_y\bigl(\vb{g}[u,r]\bigr)\circ \Phi_{\vb{g}^{-1}}\colon [u,gv]\to \Psi_{\vb{g}}\mu_y\bigl(\vb{g}[u,r]\bigr)[u,v]=\Psi_{\vb{g}}[u,r(v)]=[ug,r(v)]\;,\]
and so $\mu_y$ is equivariant.
\end{proof}
There are some special cases of Proposition~\ref{prop:MuIsIsomorphism} worth mentioning.
\begin{itemize}
\item $Y=G$ with the $G$ action induced by left multiplication. Then $\pi^*P\times_G G$ can be identified with $\pi^*P$. Thus, we see that the gauge action on sections of $\pi^*P$ is induced by the standard gauge action on $\vb{C^\infty(V,G)}$ via the isomorphisms $\mu_U$.
\item $Y=\R$. Then Proposition~\ref{prop:MuIsIsomorphism} allows us to identify smooth functions on $X$ with sections of $\vb{C^\infty(V)}$.
\item $Y=T\oplus V$. Then $\pi^*P\times_G Y$ is the tangent bundle $TX$, and vector fields on $X$ can be identified with sections of $\vb{C^\infty(V)\otimes(T\oplus V)}$.
\item $Y=\Lambda^*(T\oplus V)$. Then
$C^\infty(V,Y)=C^\infty(V,\Lambda^*V)\otimes\Lambda^*T$, which can be identified with the space $\Omega(V,\Lambda^*T)$ of $\Lambda^*T$\nobreakdash-valued differential forms on $V$ in a $G$-equivariant way. Notice that we are thinking of $V$ as a manifold when writing $\Omega(V)$, but we shall refrain from calling $\Lambda(V)$ the bundle of differential forms over $V$, which could lead to some confusion.
Summing up, we have constructed an equivariant isomorphism
\[\mu_M: \Gamma(\vb{\Omega(V\otimes\Lambda^*T)})\to\Omega(X)\;.\]
\end{itemize}
\begin{corollary}
\label{cor:InvariantForms} Under the identification $\mu_M$, the action of $\mathcal{G}$ on $\Omega(X)$ is
given by
\[[u,g][u,\alpha\otimes\beta]=[u,g\alpha\otimes g\beta]\;.\]
In particular, the space of invariant forms is
\[\Omega(X)^{\mathcal{G}}=C^\infty(M)\otimes_\R\mu_M\Omega(V,\Lambda^*T)^G\;.\]
\end{corollary}
\begin{proof}
By construction,
\[[u,g][u,\alpha\otimes\beta]=[ug,\alpha\otimes \beta]=[u,g\alpha\otimes g\beta]\;,\]
proving the first part of the theorem. As a consequence, the space of invariant forms is the image under $\mu_M$ of
\[\Gamma(P\times_G (\Omega(V,\Lambda^*T)^G))\;.\]
On the other hand, $\vb{\Omega(V,\Lambda^*T)^G}$ is a trivial bundle, because by construction the structure group $G$
acts trivially on its fibre. Therefore
\[ \Gamma(\vb{\Omega(V,\Lambda^*T)^G})=C^\infty(M)\otimes_\R\Omega(V,\Lambda^*T)^G\;,\] concluding the proof.
\end{proof}
The space of invariant forms we have found does not have the desirable property of being closed under $d$.
Indeed, it contains all functions on $M$, but it does not contain all exact forms on $M$. In the next section
we shall see that under certain hypotheses on the connection form $\omega$, the space of parallel invariant forms is
$d$-closed.

\section{$d$-closedness}
\label{sec:Closure} We have seen that the algebra of invariant forms is not closed under $d$. Having fixed a connection form $\omega$ on $P$, we can work around this
problem by discarding the $C^\infty(M)$ factor, and consider
\[\mu_M\Omega(V,\Lambda^*T)^G\subset \Omega(X)^{\mathcal{G}}\;.\]
It turns out that this subalgebra of $\Omega(X)$ may be closed under $d$, but only under very restrictive conditions on the connection, determining which is the purpose of this section.

First, we characterize the subalgebra under consideration as follows:
\begin{lemma}
\label{lemma:InvariantParallel} A section of $\vb{\Omega(V,\Lambda^*T)}$ lies in $\Omega(V,\Lambda^*T)^G$ if and
only if it is both $\mathcal{G}$-invariant and parallel.
\end{lemma}
\begin{proof}
A section of $\vb{\Omega(V,\Lambda^*T)}$ can be viewed as a $G$-equivariant map
\[\alpha:P\to \Omega(V,\Lambda^*T)\;;\]
the section is parallel if and only if
\[d\alpha+\omega\alpha=0\;.\]
This equation implies that $\alpha$ takes values in a single $G$-orbit; if the corresponding form is $\mathcal{G}$-invariant, by Corollary~\ref{cor:InvariantForms} the map $\alpha$ takes values in $\Omega(V,\Lambda^*T)^G$, and so it is constant. The converse is proved in a similar way.
\end{proof}

The second ingredient will be a formula that gives the action of $d$ on $\Omega(X)$ in terms of the bundle $\vb{\Omega(V,\Lambda^*T)}$. Some language is needed to write this formula.
We define the covariant derivative
\[\nabla\colon\Gamma(\vb{\Omega(V,\Lambda^*T)})\to \Omega^1(M,\vb{\Omega(V,\Lambda^*T)})\] in the usual way.
The natural map
\[\Lambda^pT\otimes\Lambda^qT\ni \alpha\otimes\beta\to \alpha\wedge\beta\in\Lambda^{p+q}T\]
induces a map $c$ through the composition
\[\xymatrix{
\Omega^p(M,\Omega(V,\Lambda^*T))\ar@/_7mm/[rr]^c\ar[r]^{\cong\quad}& \Gamma(\vb{\Lambda^pT\otimes\Omega(V,\Lambda^*T)})\ar[r]&\Gamma(\vb{\Omega(V,\Lambda^*T)})}\]
We also have the exterior covariant differential
\[D\colon\Omega^p\bigl(M,\vb{\Omega(V,\Lambda^*T)}\bigr)\to \Omega^{p+1}\bigl(M,\vb{\Omega(V,\Lambda^*T)}\bigr),\]
which coincides with $\nabla$ for $p=0$ and satisfies the Leibnitz rule
\[D(\alpha\otimes h)= d\alpha\otimes h +(-1)^{\deg\alpha}\alpha\wedge\nabla h\;.\]
Finally, the standard operator
$d$ on $\Omega(V,\Lambda^*T)$
is $G$-equivariant; as such, it induces bundle maps
\[\vb{d}\colon \vb{\Omega^p(V,\Lambda^*T)}\to\vb{\Omega^{p+1}(V,\Lambda^*T)}\;.\]
\begin{lemma}
\label{lemma:dOnInvariantFunctions} The map $d:C^\infty(X)\to\Omega^1(X)$ satisfies
 \begin{equation}
 \label{eqn:dOnInvariantFunctions}
d(\mu_Mh)=\mu_M\vb{d}h+(\mu_M c\nabla h),
\end{equation}
 where $h$ is a section of $\vb{C^\infty(V)}$.
\end{lemma}
\begin{proof}
It is sufficient to prove the statement at a point $[u,v]$ of $X$; let $s$ be a local section of $P$ with
$u=s(y)$ for some $y$ in $M$, and write $h=[s,f]$. We must evaluate both sides of
(\ref{eqn:dOnInvariantFunctions}) on a vector \[[(u,v),w]\in T_{[u,v]}X\;,\quad w\in T\oplus V\;.\]
Suppose first
that $w$ lies in $V$; then $\nabla h$ gives no contribution. Let $\sigma$ be a curve on $V$ with
\[\sigma(0)=v\;,\quad \sigma'(0)=w\;;\] then $\tilde\sigma=[u,\sigma]$ is a curve on $X$ satisfying
$\partial_t\tilde\sigma(0)=[(u,v),w]$.
Clearly, we have
\[d(\mu_Mh)\bigl([(u,v),w]\bigr)=\partial_t(\tilde\sigma^*\mu_Mh)(0)=df_v(w);\]
since $\mu_M\vb{d}h=\mu_M(s,df)$, (\ref{eqn:dOnInvariantFunctions}) is satisfied when evaluated on a vertical $w$.

Now let $w$ lie in $T$; then $\mu_M(\vb{d}h)$ gives no contribution. Take a horizontal curve $\sigma(t)$ on $P$ satisfying
\[\sigma(0)=u\;,\quad \theta_u\sigma'(0)=w\in T.\]
Then $\tilde\sigma=[\sigma,v]$ is a curve on $X$, and the isomorphism \eqref{eqn:lambda} maps
\[(u,v;\sigma'(0),0)\to\left((u,v),(w,0)\right) ,\]
meaning that $\partial_t\tilde\sigma(0)=[(u,v),w]$.
We can assume that $s$ agrees with $\sigma$, in the sense that $s_{p(\sigma(t))}=\sigma(t)$. In particular, $s_{*y}([u,w])$ is horizontal.
By the general formula for the covariant derivative,
\begin{equation}
 \label{eqn:dOnInvariantFunctions-2}
\nabla_{[u,w]} [s,f]=\left[u,\partial_{\,[u,w]} f + s^*\omega([u,w])f\right]=\left[u,\partial_{\,[u,w]}
f\right]\;,
\end{equation}
where $\partial_{\,[u,w]} f$ is the derivative along $[u,w]$ of $f\colon M\to C^\infty(V)$. Thus
\[\mu_M c\nabla h ([(u,v),w])=\partial_{\,[u,w]}f(v)=\partial_t(\tilde\sigma^*\mu_M h)(0)=d(\mu_Mh)([(u,v),w]),\]
where the second equality follows from $(\mu_Mh)(\sigma(t),v)=f_{p(\sigma(t))}(v)$.
\end{proof}
Recall that the torsion $\Theta\in\Omega^2(P,T)$ is the horizontal part of $d\theta$, which can also be seen as a $2$-form in $\Omega^2(M,TM)$ satisfying
\[\Theta(X,Y)=\frac12\left(\nabla_X Y-\nabla_Y X-[X,Y]\right).\]
More generally, if $S$ is a $p$-form taking values in $TM$ and $\alpha$ is a $q$-form on $M$, we define the contraction
\begin{multline*}
(S\hook_M \alpha)(X_1,\dotsc,X_{p+q-1})=\\
=\frac{q}{(p+q-1)!}\sum_{\sigma} (-1)^{\sigma}\alpha(S(X_{\sigma_1},\dotsc,X_{\sigma_p}),X_{\sigma_{p+1}},\dotsc,X_{\sigma_{p+q-1}})
\end{multline*}
where $\sigma$ runs over permutations of $p+q-1$ elements.
\begin{lemma}
\label{lemma:cnabla}
 Let $h$ be an element of $\Omega\bigl(M,\vb{\Omega(V,\Lambda^*T)}\bigr)$. Then
\begin{equation*}
\nabla_Y c h = c\nabla_Y h,\quad c\nabla c h = cD h - c(\Theta\hook_M h)\;.
\end{equation*}
\end{lemma}
\begin{proof}
We work on the principal bundle $P$, adopting the language and methods of \cite{KobayashiNomizu}. We establish the following notation: if $Y$ is a vector field on $M$, $Y^*$ is its horizontal lift. To a section  $h$ of an associated bundle, say $P\times_G W$, we associate the corresponding equivariant map $\tm{h}\colon P\to W$. More generally, to a form $h$ in $\Omega(M,\vb{W})$, we associate a form $\tm h\in\Omega(P,W)$  characterized by
\[\tm h(X_1^*,\dotsc,X_k^*)=\tm{h(X_1,\dotsc, X_k)},\quad  Z\hook\tm h =0\text{ for $Z$ vertical}.\]
The covariant derivative of a section $h$ or a vector field $X$ satisfies
\begin{equation}
\label{eqn:CovariantDerivativeOfTensorialForm}
\tm{\nabla_Y h}=Y^*\tm h,\quad \theta(\nabla_Y X)^*=Y^*\theta(X^*).
\end{equation}

At every point of $P$, the tautological form $\theta$ defines an isomorphism of the horizontal distribution with $T$. This induces a global basis of horizontal vector fields $E_1,\dotsc, E_n$. Let $e_i=\theta(E_i)$ be the corresponding basis of $T$. If $h$ is a $\vb{\Omega(V,\Lambda^*T)}$\nobreakdash-valued form on $M$, as in the statement, then
\[\tm{c h}=\sum_{1\leq i_1,\dotsc,i_k\leq n} e^{i_1,\dotsc,i_k} \wedge \tm{h}(E_{i_1},\dotsc, E_{i_k}).\]
Taking the covariant derivative with respect to $Y$,
\begin{multline*}
\tm{\nabla_Y ch}= \sum_{\substack{1\leq i_1,\dotsc,i_k\leq n}} e^{i_1,\dotsc,i_k}\wedge Y^*(\tm h(E_{i_1},\dotsc, E_{i_k}))\\
=\sum_{\substack{1\leq i_1,\dotsc,i_k\leq n}}  e^{i_1,\dotsc,i_k}\wedge (\Lie_{Y^*}\tm h)(E_{i_1},\dotsc, E_{i_k})\\
+\sum_{\substack{1\leq i_1,\dotsc,i_k\leq n\\l=1,\dotsc,k}} e^{i_1,\dotsc,i_k}\wedge (-1)^{l+1} \tm h([Y^*, E_{i_l}],E_{i_1},\dotsc,\widehat{E_{i_l}},\dotsc, E_{i_k})\,.
\end{multline*}
By the definition of $\tm{\Theta}=d\theta$ and \eqref{eqn:CovariantDerivativeOfTensorialForm},
\[\theta([Y^*, E_{i_l}])=-2\tm{\Theta}(Y^*,E_{i_l})-\tm{\nabla Y}(E_{i_l}).\]
In general, if $\alpha$ is a $TM$-valued one-form, then
\begin{multline*}
\tm{c(\alpha\hook_M h)}= \sum_{\substack{1\leq i_1,\dotsc,i_k\leq n}} e^{i_1,\dotsc,i_k}\wedge
\frac{1}{(k-1)!}\sum_{\sigma}(-1)^{\sigma} h(\alpha(E_{\sigma_1}),E_{\sigma_2},\dotsc,E_{\sigma_k})\\
=k\sum_{\sigma} e^{\sigma_1,\dotsc,\sigma_k}\wedge
 h(\alpha(E_{\sigma_1}),E_{\sigma_2},\dotsc,E_{\sigma_k})\;.
\end{multline*}
Therefore,
\begin{equation}
\label{eqn:nabla_Ych}
\tm{\nabla_Y ch}=c(\Lie_{Y^*}\tm h)-\tm{c((Y\hook \Theta)\hook_M h)}-\tm{c((\nabla Y)\hook_M h)}\;.
\end{equation}

Similarly, from
\begin{multline*}
\tm{\nabla_Y h}(X_1^*,\dotsc,X_k^*)={Y^*} (\tm h(X_1^*,\dotsc, X_{k}^*))\\-\sum_{l=1}^k (-1)^{l+1} \tm{h}(Y^*\theta(X_l^*),X_1^*,\dotsc,\widehat{X_l^*},\dotsc,X_k^*)
\end{multline*}
we compute
\[\tm{\nabla_Y h}=(\Lie_{Y^*} \tm{h})-\tm{(\nabla Y)\hook_M h}-\tm{(Y\hook \Theta)\hook_M h}\;;\]
comparing with \eqref{eqn:nabla_Ych}, the first part of the statement follows. The second part is proved in the same way.
\end{proof}

We need to introduce two more contractions beside $\hook_M$. The first is $\hook_T$, induced by the usual interior product
\[\hook:T\otimes \Omega(V,\Lambda^kT)\to\Omega(V,\Lambda^{k-1}T)\;;\]
more precisely, if $\alpha$ is a $TM$-valued $p$-form and $h$ is a section of $\vb{\Omega(V,\Lambda^*T)}$,
\[(\alpha\hook_T h)_u(X_1,\dotsc,X_p)=\alpha_u(X_1,\dotsc,X_p)\hook h_u\;.\]
Likewise, the $G$-equivariant map $\eta\colon \lie{g}\otimes\Omega(V,\Lambda^*T)\to\Omega(V,\Lambda^*T)$,
\[ \eta(A,h)(v)=A(v)\hook h(v), \quad A\in\lie{g},\; h\in\Omega(V,\Lambda^*T), \]
induces a contraction, also to be denoted by $\eta$. Explicitly, if $\alpha$ is a $\vb{\lie{g}}$\nobreakdash-valued $p$-form and $h$ is a section of $\vb{\Omega(V,\Lambda^*T)}$,
\[\eta(\alpha,h)_u(X_1,\dotsc,X_p)=\eta(\alpha_u(X_1,\dotsc,X_p),h_u).\]
We can finally prove:
\begin{lemma}
\label{lemma:Formulazza}
The map $d\colon\Omega(X)\to\Omega(X)$ satisfies
\begin{equation}
\label{eqn:Formulazza}
d(\mu_M h)=\mu_M\vb{d}h+\mu_Mc\nabla h+\frac12\mu_Mc(\Theta\hook_T h)+\frac12\mu_M c\eta(R,h)\,.
\end{equation}
\end{lemma}
\begin{proof}
We proceed by induction on the  degree of $\mu_Mh$. If the degree is zero, \eqref{eqn:Formulazza} is the content of Lemma~\ref{lemma:dOnInvariantFunctions}. Otherwise, it suffices to verify it  on a covering of $X$ consisting of subsets $X_A=\pi^{-1}(A)$, with $A$ open in $M$. We can assume that $X_A$ is diffeomorphic to some $\R^n$, and therefore assume that
\[\mu_M h=d\mu_Ma\wedge\mu_M\beta\;,\]
with  $a$ a section of $\vb{C^\infty(V)}$ and $\beta$ a  section of $\vb{\Omega(V,\Lambda^*T)}$. Indeed, as $X_A\cong\R^n$, the general element $h$ will be a linear combination with constant coefficients of elements of this type.

By induction, the left-hand side of \eqref{eqn:Formulazza} is
\begin{equation}
\label{eqn:Formulazza:LHS}
d(d\mu_Ma\wedge\mu_M\beta)=-d\mu_Ma\wedge \mu_M\left(\vb{d}\beta+ c\nabla\beta+\frac12\Theta\hook_T\beta+\frac12c\eta(R,\beta)\right)\;.
\end{equation}
On the other hand, by Lemma~\ref{lemma:dOnInvariantFunctions},
\[d\mu_Ma\wedge\mu_M\beta=\mu_M(\vb{d}a\wedge \beta+c\nabla a\wedge\beta)\;,\]
so the right hand side of \eqref{eqn:Formulazza} is (the image under $\mu_M$ of)
\begin{multline*}
\vb{d}\bigl(\vb{d}a\wedge\beta+c\nabla a\wedge\beta\bigr)+c\nabla(\vb{d}a\wedge\beta+c\nabla a\wedge\beta)\\+
\frac12\Theta\hook_T(\vb{d}a\wedge \beta+c\nabla a\wedge\beta)+\frac12c\eta(R,\vb{d}a\wedge \beta+c\nabla a\wedge\beta)
\end{multline*}
which, since $\vb{d}$ is a differential operator and it commutes with $\nabla$, differs from \eqref{eqn:Formulazza:LHS} by
\begin{equation}
\label{eqn:Difference}
c\nabla c\nabla a\wedge\beta+\frac12\Theta\hook_T(c\nabla a)\wedge\beta+\frac12 c\eta(R,\vb{d}a)\wedge\beta\;.
\end{equation}
Now, by Lemma~\ref{lemma:cnabla}
\[c\nabla c\nabla a(X,Y)=\frac12(R(X,Y)a-\nabla_{\Theta(X,Y)}a)\]
and
\[\eta(R,\vb{d}a)(X,Y)_v=\vb{d}a(R(X,Y)_v)=(-R(X,Y)a)_v,\]
so that \eqref{eqn:Difference} is zero.
\end{proof}
Recall from the beginning of this section that we are interested in sections of $\vb{\Omega(V,\Lambda^*T)}$ corresponding to elements of
$\Omega(V,\Lambda^*T)^G$. More generally, given a $G$-space $Y$, we say that the space of \dfn{constant} sections of $\vb{Y}$ is
\[Y^G\subset \Gamma(\vb{Y})\;.\]
We can now prove the main result of this section, that gives a sufficient condition on the connection $\omega$ for the algebra of constant forms to be closed.

\begin{theorem}
\label{thm:ParallelTorsion}
If both torsion and curvature are constant, the space $\Omega(V,\Lambda^*T)^G$ of constant forms is a
differential graded subalgebra of $\Omega(X)$.
\end{theorem}
\begin{proof}
We must prove that for every constant form $h$, $d\mu_h$ is also constant. By Lemma~\ref{lemma:InvariantParallel}, $h$ is $\mathcal{G}$-equivariant and parallel; by \eqref{eqn:Formulazza}, we must prove the same holds for
\[\vb{d}h+\frac12c(\Theta\hook_T h)+\frac12 c\eta(R,h)\;.\]
Since $\vb{d}$ and the contractions are equivariant, equivariance is straightforward. Concerning parallelism,
Lemma~\ref{lemma:cnabla} gives immediately
\[
\nabla_Y\left(\vb{d}h+\frac12c(\Theta\hook_T h)+\frac12 c\eta(R,h)\right)=\frac12c(\nabla_Y\Theta)\hook_T h+
\frac12 c\eta(\nabla_Y R,h))=0 \,.\qedhere
\]
\end{proof}
\begin{remark}
By way of converse, we now show that if $d$ preserves the space of constant forms, under the only additional hypothesis that the action of $\lie{g}$ on $V$ is effective, then the curvature $R$ must be constant.
Indeed suppose that $\beta$ is any non-zero form in $\Omega^2(M,\vb{\lie{g}})$.
If $h$ is the constant section of $\vb{\Omega(V,\Lambda^*T)}$ corresponding to the volume form on $V$, we claim that
\begin{equation}
\label{eqn:thm:ParallelTorsion:1}
c\eta(\beta,h)\text{ is not identically zero.}
\end{equation}
Indeed, since $\lie{g}$ acts effectively on $V$, there exist $X_1$, $X_2$, $u$ and $v$ such that
\[\beta_y(X_1,X_2)[u,v]\neq 0,\quad X_1,X_2\in T_yM, u\in P_y, v\in V\;.\]
 Then
\[\eta(\beta,h)(X_1,X_2)_{[u,v]}=\beta(X_1,X_2)[u,v]\hook h([u,v])\neq 0.\]
So there exist sections $Y_1,\dotsc,Y_k$ of $\vb{V}$ such that
\[c\eta(\beta,h)_y(X_1,X_2,Y_1,\dotsc,Y_k)=h(\beta_y(X_1,X_2)[u,v],Y_1,\dotsc,Y_k)\neq 0,\]
that implies \eqref{eqn:thm:ParallelTorsion:1}.

Now observe that by construction $\vb{d}h$ and $\Theta\hook_T h$ are zero, so that $d\mu_M h$ can only be parallel if $c\eta(\nabla_Y R, h)$ is zero. By \eqref{eqn:thm:ParallelTorsion:1}, this implies that $\nabla_Y R=0$. The same argument applied to $\beta=R-\vb{g}R$ proves that $R$ must be equivariant.
\end{remark}

An analogous statement for the torsion seems more complicated to prove: indeed, we cannot repeat the argument using the volume form on the base, because this is automatically closed. However, assume that there is an invariant element
\[s\otimes f\in S^*V\otimes \Hom(T,\Lambda^kV)\;,\]
with $f$ injective, as is the case in many interesting examples\ (see e.g. Section~\ref{sec:Examples}).
Let $t_1,\dotsc,t_n$ be an orthonormal basis of $T$; then $s\otimes f$ induces an invariant element
\[V \ni v\to s(v)\sum_{i=1}^n t_i\otimes f(t_i)\in \Lambda^{k+1}(T\oplus V)\;.\]
Let $h$ be the corresponding constant section, taking values in $\Omega^{k}(V,T)$. If $d\mu_Mh$ is constant, then by \eqref{eqn:Formulazza} it takes values in
\[\Omega^{k+1}(V,T)\oplus\Omega^{k}(V,\Lambda^2T)\oplus\Omega^{k-1}(V,\Lambda^3T)\;.\]
These spaces are invariant under both  $\mathcal{G}$ and $\nabla_Y$, so in particular the second component, namely
$\frac12c(\Theta\hook_T h)$, must be constant. Proceeding as in the case of $R$, we find that $\Theta$ must also be constant.

\smallskip
Recall that our final goal is to provide a method to construct explicit special geometries. If the algebra of constant forms is  $d$-closed, one can choose a finite basis of generators  and compute $d$ in terms of that basis (see Section~\ref{sec:Algorithm}).  Then, searching for a special geometry whose defining forms are invariant becomes a matter of looking for differential forms of the correct algebraic type, on which the linear operator $d$ satisfies the required conditions (see Section~\ref{sec:Examples}).  Otherwise, one can still hope to find single closed constant elements applying Lemma~\ref{lemma:Formulazza}, though we shall not pursue this method  in this paper.

\section{Invariant forms on a vector space}
We are now confronted with the purely algebraic problem of determining the
space of $G$-invariant $\Lambda^*T$-valued differential forms on $V$. In this section we give an effective criterion to determine whether a space of invariant forms is complete. The only assumption will be that $G$ is compact and connected.

Throughout this section, we shall identify $\Omega(V,\Lambda^*T)^G$ with
\[\mathcal{F}=\{\alpha:V\to \Lambda\mid \alpha \text{ is smooth and $G$-equivariant}\}\;, \quad \Lambda=\Lambda^*T\otimes\Lambda^*V\;.\]
We view $\mathcal{F}$ as a $C^\infty(V)^G$-module. For every subspace $W$ of $\mathcal{F}$ and $v$ in $V$, we set
\[W_v=\{\alpha(v)\mid \alpha\in W\}\;.\]
\begin{lemma}
\label{lemma:CompleteDictionary} Let $W\subset W'$ be submodules of $\mathcal{F}$, such that
\[W_v=W'_v\;,\quad\forall v\in V\;.\]
Assume further that $W$ is finitely generated. Then $W=W'$.
\end{lemma}
\begin{proof}
Suppose $W$ is finitely generated. We shall construct $w_1,\dotsc,w_k$ in $W$ such that
\begin{equation}
\label{eqn:CompleteDictionary:Generators}
\{(w_i)(v)\mid i=1,\dotsc,k,(w_i)(v)\neq 0\} \text{ is a basis of } W_v \quad\forall v\in V\;,
\end{equation}
and prove that $w_1,\dotsc,w_k$ generate $W'$. Let $w_1,\dotsc,w_k$ be any
set of generators of $W$. Consider the  set \[U=\left\{v\in V\mid w_k(v) \text{ is linearly dependent on }
w_1(v),\dotsc, w_{k-1}(v)\right\}\;,\] and let $f$ be a non-negative smooth function on $V$ with $f^{-1}(0)=U$. For
instance, we can take $f$ to be the square distance between $w_k(v)$ and the vector space spanned by
$w_1(v),\dotsc, w_{k-1}(v)$. Now define
$$h(x)=\int_G f(gx) \mu_G\;,$$ where $\mu_G$ is the Haar measure on $G$. Clearly, $h$ is a $G$-invariant function with
$h^{-1}(0)=U$. Thus, we have obtained elements
$$w_1,\dotsc,w_{k-1},h\,w_k\in W$$ such that at each $v$, the last element is either zero or linearly
independent of $w_1,\dotsc,w_{k-1}$. Iterating this procedure for $w_{k-1},\dotsc,w_2$, we
find elements satisfying \eqref{eqn:CompleteDictionary:Generators}.

Elements of $W'$ can be viewed as sections of a trivial bundle $\pi\colon V\times \Lambda\to V$. The $w_i$ define a degenerate metric on this bundle by
\[w_1\otimes w_1+\dotsb+w_k\otimes w_k\;.\]
Since $W_v=W'_v$, for every $\alpha$
in $W'$
\[\alpha(v)=\sum_{i=1}^k \langle \alpha(v),w_i(v)\rangle w_i(v)\;,\] i.e. $\alpha=\sum \langle\alpha,w_i\rangle w_i$.
By construction, $\langle\alpha,w_i\rangle$ are smooth invariant functions; so, $\alpha$ is linearly dependent on the $w_i$.
\end{proof}
We write $H\leq G$ if $H$ is a closed subgroup of $G$.
We shall need the following result from \cite{Bredon}:
\begin{theorem}[Tietze-Gleason-Bredon]\label{thm:TietzeGleasonBredon}
Let $M$ and $N$ be manifolds on which a compact Lie group $G$ acts, and let $A$ be a closed invariant subspace of $M$.
Then every smooth equivariant map $\phi:A\to N$ can be extended to a smooth equivariant map $\overline{\phi}:M\to N$.
\end{theorem}

We can now prove the main theorem of this section:
\begin{theorem}
\label{thm:CompleteDictionary} $\mathcal{F}$ is finitely generated. Let $W$ be a finitely generated submodule of
$\mathcal{F}$; then $W$ coincides with all of $\mathcal{F}$ if and only if
\begin{equation}
\label{eqn:CompleteDictionary} W_v=\Lambda^{\Stab v}\; \forall v\in V\;.
\end{equation}
\end{theorem}
\begin{proof}
Let $W$ be a finitely generated submodule of $\mathcal{F}$ satisfying (\ref{eqn:CompleteDictionary}). Let
$\alpha\colon V\to \Lambda$ be invariant; clearly, $\alpha(v)$ is fixed by $\Stab(v)$. So,
\begin{equation}
\label{eqn:WFS}
W_v\subseteq\mathcal{F}_v\subseteq \Lambda^{\Stab(v)}\;,
\end{equation}
and by hypothesis equality holds. Then  $W=\mathcal{F}$ by Lemma \ref{lemma:CompleteDictionary}, proving the ``if'' part of the statement.

We claim that there always exists $W$ as in the hypothesis. Then $\mathcal{F}=W$ is finitely generated, and equality holds in \eqref{eqn:WFS}, giving the ``only if'' part of the statement.

To prove the claim,
choose a group $H$ which is the stabilizer of a point of $V$, and let $s_1,\dotsc,s_k$ be a basis of $\Lambda^{H}$. For each
$s_i$, we shall construct an equivariant map
$$\tilde{s}_i:V\to \Lambda$$ such that $\tilde{s}_i(v)$ is a non-zero multiple of $s_i$ at each $v$ with stabilizer $H$. By
equivariance, it will follow that $\tilde{s}_1(v),\dotsc,\tilde{s}_k(v)$ is a basis of $\Lambda^{\Stab v}$ for all $v$ with
stabilizer conjugate to $H$. It is well known that an orthogonal representation of a compact group has a finite number of orbit types (see \cite{Bredon}); in other words, there are only a finite number of stabilizers up to conjugation. It is therefore sufficient to consider a finite number of
$H$; thus, we obtain a finite collection of equivariant maps from $V$ to $\Lambda$, which generate a module $W$ with the
required properties. It only remains to construct the $\tilde{s}_i$.

Observe first that since $G$ is compact, the number of $K\leq G$ which have the same Lie algebra as $H$ is finite. For
every such $K$ containing $H$, choose an element $g_K\in K\setminus H$. At each point $v$ of $V^H$, consider the map
$$\lie{g}/\lie{h} \to \End(V)$$ induced by the infinitesimal action of $G$; this map is injective if and only if $\Stab
v$ has the same Lie algebra as $H$. Therefore we can define a smooth function $f$ on $V^H$ satisfying
\begin{equation*}
\begin{cases} f(v)>0 &\text{ if } \Stab v \text{ has the same Lie algebra as }H\;,\\
                f(v)=0 &\text{otherwise.}
\end{cases}
\end{equation*}
Now consider the function on $V^H$
\[\tilde{f}(v)=f(v)\prod_{\substack{\lie{k}=\lie{h}\\H\leq K\leq G}}\norm{g_K(v)-v}^2\;;\]
by construction, $\tilde{f}^{-1}(0)$ is the subset of $V^H$ where the stabilizer is bigger than $H$. So, we can define a map
\[ V^H\ni v\xrightarrow{\tilde{s}_i} \tilde{f}(v)s_i\in \Lambda\;,\]
which extends to an equivariant map on $GV^H$. Since $GV^H$ is closed, we can apply Theorem
\ref{thm:TietzeGleasonBredon} to extend $\tilde{s}_i$ to an equivariant map defined on all of $V$.
\end{proof}
\begin{remark}
When checking the criterion of Theorem \ref{thm:CompleteDictionary}, it may be useful to separate forms by degree. More
precisely, we say that a form in $\mathcal{F}$ has \dfn{bidegree} $(p,q)$ if it takes values in
\[\Lambda_{p,q}=\Lambda^pT\otimes\Lambda^qV\;;\] we denote
$\mathcal{F}_{p,q}$ the space of forms of bidegree $(p,q)$. Since $\Lambda=\bigoplus \Lambda_{p,q}$ is a $G$-invariant decomposition, we have thus made $\mathcal{F}=\bigoplus \mathcal{F}_{p,q}$ into a bigraded vector
space. Theorem \ref{thm:CompleteDictionary} still holds if one replaces
(\ref{eqn:CompleteDictionary}) with
\[(W\cap\mathcal{F}_{p,q})_v=(\Lambda_{p,q})^{\Stab v} \quad \forall v\in V\;.\]
\end{remark}

We now apply Theorem~\ref{thm:CompleteDictionary} to show that the dictionary of  Section~\ref{sec:TS2} is complete. In this case, $G=\LieG{U}(1)$ and $V=\R^2=T$. We can reinterpret the functions $a_1$, $a_2$ on $T^*S^2$ as  the standard coordinates on $V=\R^2$, and, accordingly, identify the $b_i$ with $da_i$. Moreover, we declare $\beta_1$, $\beta_2$ to be the standard basis of $T^*=(\R^2)^*$. We can now state and prove:
\begin{corollary}
The $C^\infty(V)^G$-algebra generated by \eqref{eqn:contractions}, \eqref{eqn:contractions2} and \eqref{eqn:contractions3} coincides with the space of invariant forms $\Omega(V,\Lambda^*T)^G$.
\end{corollary}
\begin{proof}
 The stabilizer for the action of $G$ on $V$ is the whole $G$ at $v=0$, and the identity at other points. By a straightforward computation,
\[\dim (\Lambda_{p,q})^{G}=\begin{cases} 1&p=0,2, \;q=0,2,\\ 2&p=1=q,\\ 0&\text{otherwise.}\end{cases}\]
Let $W$ be the $C^\infty(V)^G$-algebra generated by \eqref{eqn:contractions}, \eqref{eqn:contractions2} and \eqref{eqn:contractions3}.
Then $W_0$ is
\[W_0=\left\langle 1,b_{12},\beta_{12},b_1\beta_1+b_2\beta_2,b_1\beta_2-b_2\beta_1, b_{12}\wedge\beta_{12}\right\rangle,\]
which has dimension six; so, \eqref{eqn:CompleteDictionary} is satisfied at $v=0$. Evaluating at the point $v$ with coordinates $a_1=1$, $a_2=0$, we obtain
\[W_v\supset\langle b_1, \beta_1, b_2,\beta_2\rangle\;.\]
Thus, $W_v=\Lambda$, and Theorem~\ref{thm:CompleteDictionary} is satisfied.
\end{proof}

\section{Dictionaries of forms}
\label{sec:Dictionary}
Theorem~\ref{thm:ParallelTorsion} and the following remarks show that the natural setting for our construction is that of constant torsion and curvature. This condition implies that any two points have affine isomorphic neighbourhoods (see \cite{KobayashiNomizu}, Chap. VI);  this is very close to requiring that  the group of affine transformations act transitively. In this section we shall assume that this stronger condition holds, considering, that is to say,  homogeneous spaces, and we illustrate the method sketched in Section~\ref{sec:TS2} to produce a list of elements of $\Omega(V,\Lambda^*T)^G$. This list is called a ``dictionary'', because its elements are generated from a much smaller list of objects in the same way that the words appearing in a dictionary are obtained from the letters of the alphabet. Consistently with the metaphor (since different words may share the same meaning), we shall think of our dictionary not so much as a list of forms, but as an abstract set that maps non-injectively to $\Omega(V,\Lambda^*T)^G$; this distinction will play a significant r\^ole in Section~\ref{sec:Algorithm}.

Consider a homogeneous space
\[M=H/G\;;\]  we view $H$ as a principal bundle over $M$ with fibre $G$, on which we fix a $G$\nobreakdash-invariant splitting
\begin{equation}
\label{eqn:hSplits}
\lie{h}=T\oplus\lie{g}\;.
\end{equation}
The projection on $\lie{g}$, extended left-invariantly to $H$, gives a connection form $\omega$; similarly, the projection on $T$ defines a tautological form $\theta$ that makes $H$ into an ``ineffective $G$-structure'': in other words,
\[TM=H\times_G T,\]
even though $G$ might not act effectively on $T$.

As usual, we fix a $G$-module $V$, setting $X=H\times_G V$, and we are interested in the algebra $\Gamma(H/G,\vb{\Omega(V,\Lambda^*T)^G})$ of constant forms in the sense of Section~\ref{sec:Closure}. Of course, there is a global action of $H$ on $X$, and it would also be natural to consider the space $\Omega(X)^H$ of forms that are invariant in this sense. In fact, the two spaces coincide:
\begin{proposition}
\label{prop:ConstantIsInvariant}
The space of constant forms on $X=H\times_G V$ coincides with the space of $H$-invariant forms, i.e.
\[\mu_M(\Gamma(M,\vb{\Omega(V,\Lambda^*T)^G}))=\Omega(X)^H.\]
\end{proposition}
\begin{proof}
We can view sections of the associated bundle $\vb{\Omega(V,\Lambda^*T)}$ as $G$-equivariant maps
\[\alpha\colon H\to\Omega(V,\Lambda^*T),\]
where $G$ acts on $H$ by right multiplication. Such a map determines a $H$-invariant form in $\Omega(X)$ if and only if, say, $\alpha(h)=\alpha_0$ for all $h$. Then, by $G$-equivariance, $\alpha_0$ is in $\Omega(V,\Lambda^*T)^G$. Conversely, every element $\alpha_0$ of $\Omega(V,\Lambda^*T)^G$ determines a unique $\alpha$ by $\alpha(h)=\alpha_0$, which under $\mu_M$ corresponds to a $H$-invariant form on $X$.
\end{proof}
Motivated by Proposition~\ref{prop:ConstantIsInvariant}, we shall study the space of $H$-invariant forms on $X$. By construction, $X$ is a quotient of $H\times V$ by $G$, acting as \[(R_g)(h,v)=(hg,g^{-1}v).\] The pull-back to $H\times V$ of a form on $X$ is a \dfn{basic} form, meaning that
\begin{enumerate}
\item{(1)} $R_g^*\alpha=\alpha$ for all $g$ in $G$.
\item{(2)} $Y\hook\alpha=0$ for all fundamental vector fields $Y$ relative to the action of $G$.
\end{enumerate}
Further, invariant forms satisfy
\begin{enumerate}
\setcounter{enumi}{2}
\item{(3)} $L_h^*\alpha=\alpha$ for all $h$ in $H$.
\end{enumerate}
Thus, instead of working with the associated bundle $\vb{\Omega(V,\Lambda^*T)^G}$, relying on Lemma~\ref{lemma:Formulazza} to determine the action of $d$, we will be concerned with forms on $H\times V$ satisfying Conditions (1)--(3), and use the fact that $d$ commutes with the pullback.

Fix a $G$-invariant metric $\sum da_i\otimes da_i$, where $a_1,\dotsc,a_k$ are linear coordinates on $V$. Let $v_1,\dotsc, v_k$ be the dual  orthonormal basis of $V$. Choose a basis $e_i$ of $\lie{h}$ consistent with \eqref{eqn:hSplits}, namely
\[T\oplus \lie{g}=\langle e_1,\dotsc e_{\dim T}\rangle \oplus \langle e_{\dim T+1},\dotsc, e_{\dim H}\rangle\,.\]
Let $e^1,\dotsc,e^{\dim H}$ be the dual basis of $\lie{h}^*$. We shall implicitly extend the $e^i$ to left-invariant one-forms on $H$.

The forms $e^i$, $1\leq i\leq\dim T$ satisfy Conditions (2) and (3). The forms $da_i$ do not satisfy Condition (2), but they give rise to forms
\[b_i = da_i+a_j(\rho_*\omega)_{ij}\]
that do satisfy Condition (2). Here, $\rho:G\to\GL(V)$ is the representation map and $\rho_*$ is its derivative at the identity.

All forms satisfying Conditions (2) and (3) can be written as
\begin{equation}
\label{eqn:GenericHorizontalForm}
\sum_{I,J}f_{IJ}(a_1,\dotsc a_k) e^I\wedge b^J,
\end{equation}
where  $I$, $J$ are multiindices relative to the ranges $1\leq i\leq\dim T$, $1\leq j\leq k$ respectively.
However, the generic form \eqref{eqn:GenericHorizontalForm} does not satisfy condition 1.

The key observation is that whilst the $b_i$ taken individually do not satisfy Condition (1), the $V$-valued form $b=\sum_i b_i v_i$ is invariant, in the sense that
\begin{itemize}
\item[(1')] $R_g^*\alpha=g^{-1}\alpha$ for all $g$ in $G$.
\end{itemize}
The same holds for $a=\sum a_i v_i$. Moreover, one can contract $a$ and $b$ to obtain a legitimate invariant one-form
\[ab=\sum_i a_ib_i\;.\]
The rest of this section contains a generalization of  this procedure.

\smallskip
$V$-valued forms on $H\times V$ satisfying Conditions (1') and (2) can be identified with sections of $\Omega(X,\pi^*\vb{V})$.  The ``letters'' of our alphabet will be elements of the space
\[\mathcal{L}=(\Omega(V,\Lambda^*T)\otimes V)^G\subset \Omega(X,\pi^*\vb{V}),\]
consisting of elements of $\Omega(V\times H,V)$ satisfying Conditions (1'), (2), and (3).

As we mentioned, $\mathcal{L}$ always contains two canonical elements, namely
\begin{equation} \label{eqn:a}
a:V\to \Lambda^0T\otimes\Lambda^0V\otimes V,\quad a(v)=1\otimes v,
\end{equation}
and
\begin{equation} \label{eqn:b}
b\colon V\to\Lambda^0T\otimes\Lambda^1V\otimes V\quad b(v)=\sum_i v_i\otimes v_i.
\end{equation}
Moreover, suppose $\phi:V\to\Lambda^k T^*$ is a non-zero, $G$-equivariant map. Then there is a corresponding element of $\mathcal{L}$
\[\beta_\phi:V\to \Lambda^kT\otimes\Lambda^0V\otimes V,\quad \beta_\phi(v)=\sum_i \phi(v_i)\otimes v_i\;.\]
Lastly, suppose one has a  non-zero, $G$-equivariant map $\psi\colon V\otimes V\to\Lambda^k T^*$. An  element of $\mathcal{L}$ is induced, by
\[\epsilon_\psi:V\to \Lambda^kT\otimes\Lambda^0V\otimes V,\quad \epsilon_\psi(v)=\sum_i \psi(v,v_i)\otimes v_i\;.\]
Of course, this small list does not exhaust $\mathcal{L}$, but the point is that elements of $\mathcal{L}$ are generally easier to obtain than elements of $\Omega(V,\Lambda^*T)$.

 Elements of $\mathcal{L}$ can be combined by means of a contraction to obtain elements of $\Omega(V,\Lambda^*T)^G$. More precisely, if $m:V\otimes V\to \R$ is a $G$-equivariant map, then we can define a linear map
\begin{equation}
\label{eqn:contraction}
m_*\colon\mathcal{L}\otimes\mathcal{L}\to\Omega(V,\Lambda^*T)^G
\end{equation}
by  \[m_*\left(\sum_i\alpha_i\otimes v_i,\sum_j\alpha_j'\otimes v'_j\right)=\sum_{i,j}m(v_i,v'_j)\alpha_i\wedge\alpha'_j\;.\]
A contraction $m$ always exists, namely the scalar product. In this case, we write
\[\alpha\cdot\beta=m_*(\alpha,\beta).\]

One can then ask if the set of all possible contractions generates the whole algebra of invariant forms. This is in general not true, since $\Lambda^*T$ might contain invariant elements (corresponding to invariant forms on $M$) that cannot be obtained as contractions. However, the following holds:
\begin{proposition}
\label{prop:MethodWorks}
The algebra of invariant forms satisfies
\[(\Omega(V,\Lambda^*T))^G=\mathcal{L}\cdot\mathcal{L}+ (\Lambda^*T)^G\]
\end{proposition}
\begin{proof}
Let $\alpha$ be an invariant map
\[\alpha\colon V\to\Lambda^pT\otimes\Lambda^qV\;,\quad q>0.\]
Define an element of $\mathcal{L}$ by
\[\tilde\alpha \colon V\to\Lambda^pT\otimes\Lambda^{q-1}V\otimes V,\quad \tilde \alpha(v)= \sum_i (v_i\hook \alpha(v))\otimes v_i\]
Then
\[b\cdot\tilde\alpha=\sum_i v_i\wedge(v_i\hook a)=q\alpha\;;\]
thus, $\alpha$ lies in $\mathcal{L}\cdot\mathcal{L}$.

Now suppose $q=0$, so that $\alpha\colon V\to\Lambda^pT$. Its differential at $v$ is a map $d\alpha_v\colon V\to\Lambda^pT$, and
\[\alpha(v)=\alpha(0)+\int_{0}^1 d\alpha_{tv}(v) dt\;.\]
Introduce an element of $\mathcal{L}$ by
\[\tilde\alpha\colon V\to \Lambda^pT\otimes \Lambda^0V\otimes V,\quad
\tilde\alpha(v)=\sum_i\int_{0}^1 d\alpha_{tv}(v_i)\otimes v_i.\]
Since $\alpha(0)$ lies in $(\Lambda^*T)^G$ and
\[\alpha=\alpha(0)+a\cdot\tilde\alpha,\]
the proof is complete.
\end{proof}
Observe that the proof of Theorem~\ref{thm:CompleteDictionary} also applies to $\mathcal{L}$, showing that $\mathcal{L}$ is a finitely generated module over the space of invariant functions. Then Proposition \ref{prop:MethodWorks} poses a natural problem: having fixed $G$, $T$ and $V$, determine a finite,  minimal set $L\subset\mathcal{L}$ such that
\begin{equation}
\label{eqn:MinimalL}
(\Omega(V,\Lambda^*T))^G=C^\infty(V)^G(L\cdot L+ (\Lambda^*T)^G).
\end{equation}
This will be discussed in Section~\ref{sec:Algorithm}, though in a slightly different flavour, as we shall allow for more general contractions, beside the scalar product.

\smallskip
So far, we have not really used the fact that $M$ is a homogeneous space. This hypothesis does simplify things when computing the action of $d$.  Recall that $H\times V$ can be identified with the pullback bundle $\pi^*P$, and the connection form $\omega$ pulls back to a connection form $\omega^X$.
The connection form $\omega^X$ defines an operator $D^X$ on $\Omega(X,\pi^*\vb{V})$ by
\[D^X\alpha=d\alpha+\omega^X\wedge\alpha\;.\]
It is easy to see that $D^X$ preserves (1'), (2) and (3); thus, one can start with a small subset of $\mathcal{L}$, and enlarge it by applying $D^X$ and its powers (although $(D^X)^2$ acts like contraction with the curvature). One can then use  the Leibnitz rule
\[d(m_*(\alpha,\beta))=m_*(D^X\alpha,\beta)+(-1)^{\deg\alpha} m_*(\alpha,D^X\beta)\]
to compute the action of $d$. This method has the advantage that the ``canonical'' elements introduced before Proposition~\ref{prop:MethodWorks} satisfy some special properties with respect to $D^X$. Indeed, by definition
\[D^Xa=b\;.\]
Moreover if the torsion is zero, meaning that $H/G$ is a symmetric space, elements of the form $\beta_\phi$ are automatically $D^X$-closed. However, elements of the form $\epsilon_\psi$ do not share this property. Thus, depending on $V$, $T$ and $G$, it might be the case that the most direct method to compute the action of $d$ is pulling back the forms to $H\times V$, and compute the action of $d$ there.

\section{An algorithm}
\label{sec:Algorithm}
In this section we present an efficient algorithm to compute the dictionary of invariant forms, valid under the assumption that the group $G$ acts transitively on the sphere in $V$. This hypothesis has two practical consequences: first, $C^\infty(V)^G$ is the space of (even) functions of the radius. Second, there are only two orbit types, and in particular the set of orbit types is totally ordered. It is this latter condition that is essential to the algorithm.

The input of the algorithm is the following:
\begin{itemize}
\item A finite subset  $L\subset\mathcal{L}$.
\item A finite set  ${B}$ of contractions; more precisely, we require ${B}$ to be a disjoint union
\[{B}=\bigsqcup_{k\in\N} {B}_k,\quad{B}_k\subset\Hom(V^{\otimes^k},\R)^G.\]
\end{itemize}
Roughly speaking, the purpose of the algorithm is to check whether an analogue of \eqref{eqn:MinimalL} holds and, in that case, produce a minimal list of generators for $\Omega(V,\Lambda^*T)^G$.

Formally, we introduce the set $S$ of formal contractions of elements of $\mathcal{L}$, i.e.
\begin{equation}
\label{eqn:mathcalS}
 S=\{b(\alpha_1,\dotsc,\alpha_k)\mid b\in{B}_k, \alpha_1,\dotsc,\alpha_k\in L\}\;.
\end{equation}
The map \eqref{eqn:contraction} induces a ``translating map''
\begin{align*}
P\colon S&\to \Omega(V,\Lambda^*T)^H\\
b(\alpha_1,\dotsc,\alpha_k)&\to b_*(\alpha_1,\dotsc,\alpha_k)
\end{align*}
We can now introduce  the  associative $\R$-algebra $\mathcal{C}$ generated by $\mathcal{S}$. From the point of view of geometry, it would be natural to impose commutativity relations on $\mathcal{C}$, since its elements represent differential forms. However, our algorithm makes an essential use of a partial ordering on $\mathcal{C}$, suggesting otherwise. Also, from the point of view of  implementation there is some advantage in working with sequences, rather than products.

The translating map $P$ extends to an algebra homomorphism
\[P\colon \mathcal{C}\to\Omega(V,\Lambda^*T)^G.\]
For every $v$ in $V$, the evaluation at $v$ induces a map
\[P_v\colon \mathcal{C}\to\Lambda^*V\otimes\Lambda^*T.\]
Notice that $P$ is typically not injective. Also, $P$ cannot be surjective, because $\mathcal{C}$ was defined as the algebra generated by $S$ over $\R$, not $C^\infty(V)^G$. On the other hand, $P$ is surjective ``up to'' $C^\infty(V)^G$ if and only if $P(\mathcal{C})$ satisfies the criterion of Theorem~\ref{thm:CompleteDictionary}.  Thus, the problem of determining a minimal set of generators is addressed by the following:
\begin{algo}
\label{algo:}
The following algorithm computes a finite subset $C\subset\mathcal{C}$, of minimal cardinality, such that $\Span P(C) =P(\mathcal{C})$. \end{algo}
\begin{proof}[Description]
For every integer $l$, we define
\[S^l=\left\{ \alpha_1\dotsm\alpha_l\mid  \alpha_i\in S\right\}\subset\mathcal{C};\]
we say that elements of $S^l$ have \dfn{length} $l$.

The algorithm consists of two steps. First, one computes a set \mbox{$C_0\subset\mathcal{C}$} such that
\begin{equation}
\label{eqn:C0}\Span P_0(C_0)=P_0(\mathcal{C}), \quad\card{C_0}=\dim P_0(\mathcal{C}),\quad C_0=\bigcup_l (S^l\cap C_0).
\end{equation}
Then one fixes a non-zero $v$ in $V$, and extends $C_0$ to a set $C_v\subset\mathcal{C}$ such that
\begin{equation}
\label{eqn:Cv}
\Span  P_v(C_v)=P_v(\mathcal{C}),\quad \card{C_v}=\dim P_v(\mathcal{C}),\quad C_v=\bigcup_l (S^l\cap C_v).
\end{equation}
By Theorem~\ref{thm:CompleteDictionary}, \eqref{eqn:Cv} implies that $\Span P(C_v)=P(\mathcal{C})$ and $C_v$ is of minimal cardinality, as required.

We shall illustrate in detail the second step, as the first is completely analogous. Thus, we assume that a subset $C_0$ of $\mathcal{C}$ satisfying \eqref{eqn:C0} is given.
\begin{enumerate}
\item Fix a total ordering relation $\leq$ on $\mathcal{S}$. Extend $\leq$ to $\bigcup_l S^l$ by decreeing that
$\alpha_1\dotsm\alpha_l<\beta_1\dotsm\beta_m$ if either $l<m$ or $l=m$ and $\alpha_j<\beta_j$ for the smallest $j$ such that $\alpha_j\neq \beta_j$.
\item Compute a finite sequence of elements $\{(k_v^1)_n\}_{n>0}$, by
\begin{gather*}
(k_v^1)_n=\min\{k\in S \mid P_v(k)\notin\Span P_v(C_0\cup\{(k_v^1)_i\mid i<n\})\}.
\end{gather*}
Observe that having fixed a basis of $\Lambda^*(T\oplus V)$, this is a matter of applying an elimination scheme to a matrix with numeric entries, whose rows represent the images in $\Lambda^*(T\oplus V)$ of the elements of $S$.
\item Set $C_v^1=(S^1\cap C_0) \cup \{(k_v^1)_n\}$.
\item Starting with $l=1$, compute the subset  $C_v^{l}$ of elements of length $l$, iteratively from $C_v^{l-1}$, as follows.
\item[4.1] For each  $\alpha_1\dotsm\alpha_{l-1}\in C_v^{l-1}$,
set
\begin{multline*}
A_{\alpha_1\dotsm\alpha_{l-1}}=\biggl\{\alpha_1\dotsm\alpha_{l-1}\beta\mid \beta\in C_v^1,\alpha_{l-1}\leq\beta,
\\
\alpha_1\dotsm\widehat{\alpha_j}\dotsm\alpha_{l-1}\beta\in C_v^{l-1}\;\forall j=1,\dotsc,{l-1}\biggr\}\;.
\end{multline*}
Iterating through all elements $\alpha_1\dotsm\alpha_{l-1}$ of $C_v^{l-1}$, one obtains the set
\begin{multline*}
A_v^{l}=\bigl\{\alpha_1\dotsm\alpha_{l}\mid \alpha_j\in C_v^1,\alpha_1\dotsm\widehat{\alpha_j}\dotsm\alpha_{l}\in C_v^{l-1}\;\forall j=1,\dotsc,{l},\\
\alpha_1\leq\dotsm\leq\alpha_{l}\bigr\}.
\end{multline*}
\item[4.2] Compute a finite sequence $\{(k_v^{l})_n\}_{n>0}$ of elements of $A_v^{l}\setminus C_0$ by
\begin{multline*}
(k_v^l)_n=\min\biggl\{k\in A_v^{l}\setminus C_0 \mid \\
P_v(k)\notin\Span P_v\left(C_0\cup C_v^{1}\cup\dotsb \cup C_v^{l-1}\cup\{(k_v^l)_i\mid i<n\}\right)\biggr\}\;.
\end{multline*}
\item[4.3] Set $C_v^{l}=\{(k_v^{l})_n\}\cup (S^{l}\cap C_0)$.
\item Increment $l$ and repeat step 4, until $C_v^{l}$ is empty.
\item Set $C_v=\bigcup_l C_v^l$.
\end{enumerate}
The above steps give a construction of $C_v$ given $C_0$. On the other hand, if one replaces $C_0$ with the empty set and $v$ with $0$, the same steps give a construction of $C_0$. Let us check that the set $C_0$ obtained this way satisfies \eqref{eqn:C0}.
We show by induction on $l$ that
\begin{equation}\label{eqn:InductionStepL}
P_0(h)\in\Span \{P_0(k)\mid k\in C_0^{1}\cup\dotsb \cup C^l_0,k\leq h\},\quad  h\in S^l.
\end{equation}
The case $l=1$ is straightforward.
Define a finite sequence
\[
(h_0^l)_n=\min\bigl\{h\in (C^1_0)^{l} \mid P_0(h)\notin\Span P_0\left(C_0^{1}\cup\dotsb \cup C_0^{l-1}\cup\{(h_0^l)_i\mid i<n\}\right)\bigr\}.
\]
Here, $(C^1_0)^{l}$ is the subset of $S^l$ whose factors lie in $C^1_0$. We shall prove that
\begin{equation}
\label{eqn:InductionStepN}
(h_0^l)_n\in A_0^l\;,
\end{equation}
so that the sequences $(h_0^l)_n$ and  $(k_0^l)_n$ coincide. Since $C_0^l=\{(k_0^l)_n\}$ by definition, and
\[P_0(C_0^1)=P_0(S)\]
generates  $P_0(\mathcal{C})$ as an algebra, \eqref{eqn:InductionStepL} will follow.

Let  $h$ be the smallest element of $S^l$ such that $P_0(h)$ does not lie in
\begin{equation}
\label{eqn:HicSuntPZeri}
\Span P_0\left(C_0^{1}\cup\dotsb \cup C_0^{l-1}\cup\{(h_0^l)_i\mid i<n\}\right)\;;
\end{equation}
we must prove that $h$ lies in $A_0^l$. Let $h=\alpha_1\dotsm\alpha_l$. By minimality of $h$, we have $\alpha_1\leq\dotsb\leq\alpha_l$. Suppose by contradiction that  $\alpha_1\dotsm\widehat{\alpha_j}\dotsm\alpha_l$ is not in $C_0^{l-1}$ for some $j$. By the induction hypothesis \eqref{eqn:InductionStepL},
\[P_0(\alpha_1\dotsm\widehat{\alpha_j}\dotsm\alpha_l)=\sum_i a_iP_0(\beta_i),\; a_i\in\R, \; \beta_i\leq \alpha_1\dotsm\widehat{\alpha_j}\dotsm\alpha_l.\]
In particular,
\[P_0(h)=\sum_i \epsilon_ia_iP_0(\alpha_j\beta_i), \quad \epsilon_i=\pm1.\]
For each $i$,  the factors of $\alpha_j\beta_i$ can be reordered so as to obtain an element lesser than or equal to $h$. However, equality would imply that $\beta_i=\alpha_1\dotsm\alpha_l$, which is absurd. Since $h$ is minimal, it follows that
\eqref{eqn:HicSuntPZeri} contains all the $P_0(\alpha_j\beta_i)$, though not their linear combination $P_0(h)$, which is absurd.

Having proved \eqref{eqn:InductionStepL}, it follows that $P_0(C_0)$ spans all of $P_0(\mathcal{C})$. Notice that the termination clause of Step (5) is justified by the fact that when $C^l_0$ is empty, $A_0^{l+1}$ is also empty. This condition is bound to occur for some value of $l$ because $P_0(\mathcal{C})$ is finite dimensional. Finally, the fact that the images under $P$ of the elements of $C_0$ are linearly independent is built in the construction, so \eqref{eqn:C0} is satisfied. Equation~\eqref{eqn:Cv} is proved in the same way, bearing in mind that linear independence at $0$ implies linear independence at $v$.
\end{proof}

\begin{remark}
The fact that $S$ consists of formal contractions was not used in Algorithm~\ref{algo:}. The essential ingredient is the map $P$. Thus, one can enlarge $S$ by a set of generators for the algebra $(\Lambda^*T)^G$, extending $P$ in the obvious way. Indeed, by Proposition~\ref{prop:MethodWorks} one can always choose a finite $L$ for which, enlarging $S$ this way, Algorithm~\ref{algo:} gives a set of generators for $\Omega(V,\Lambda^*T)^G$.
\end{remark}
We conclude this section with a few words concerning the implementation of Algorithm~\ref{algo:}. Most software for differential geometry is written in the form of a package for a Computer Algebra System, see e.g. {GrTensor} \cite{GrTensor}, { GRG} and { Ricci}. This approach has the disadvantage that the underlying  language provides no support for object oriented programming, or even user-defined types. Motivated by the intrinsically polymorphic nature of the relevant data  --- the elements of $\mathcal{C}$  --- we resorted instead to an actual programming language, namely {${\rm C}\sp {++}$}, in the same vein as the library  GiNaC \cite{GiNaC}. In fact, our implementation relies heavily on GiNaC for symbolic computations, although some work was needed to introduce adequate support for differential forms.

Concerning efficiency, the algorithm was designed to minimize the size and number of the matrices appearing in Steps (2) and (4.2), by far the most computationally expensive. Empirical tests show that more na\"ive solutions lead to extremely long computation times, even for low-dimensional examples as that of Section~\ref{sec:Examples}.
\section{Calabi-Yau structures on $T\CP^2$}
In this section we construct examples of local Calabi-Yau metrics. More precisely, we consider the total space of $T\CP^2$, with the natural action of $\SU(3)$, and compute a set of generators for the space of invariant forms. Using these generators, we describe the known homogeneous $3$-Sasaki structure on the sphere subbundle $\mathbb{S}$, and the corresponding conical hyperk\"ahler structure on the complement of the zero section. As a generalization, we consider a two-parameter family of contact structures on the sphere subbundle, and, for each contact form $\alpha$, we identify the symplectic cone over the contact manifold $(\mathbb{S},\alpha)$ with the complement of the zero section in $T\CP^2$. By producing an explicit metric on $\mathbb{S}$ (which requires the list of generators) and using an extension theorem, we show that each of these symplectic manifolds has, locally, a compatible Calabi-Yau metric.
\label{sec:Examples}
\subsection{The dictionary on $\SU(3)\times_{\LieG{U}(2)}\C^2$}
We regard $\CP^2$ as the homogeneous space $\SU(3)/\LieG{U}(2)$. Accordingly, we identify its tangent space as the associated bundle
\[T\CP^2=\SU(3)\times_{\LieG{U}(2)}\C^2\;.\]
The group $\LieG{U}(2)$ acts transitively on the sphere in $V=\C^2$, and so we can apply Algorithm~\ref{algo:}.
Explicitly, we fix a basis $e_1,\dotsc,e_8$ of the Lie algebra $\su(3)$,
such that the dual basis $e^1,\dotsc,e^8$ of $\su(3)^*$  satisfies
\begin{align*}
 de^1&=-e^{23}-e^{45}+2\,e^{67},&
 de^2&=e^{13}+e^{46}-e^{57}-{\sqrt{3}}\,e^{58},\\
 de^3&=-e^{12}-e^{47}+{\sqrt{3}}\,e^{48}-e^{56}\;,&
 de^4&=e^{15}-e^{26}+e^{37}-{\sqrt{3}}\,e^{38},&\\
 de^5&=-e^{14}+e^{27}+e^{36}+{\sqrt{3}}\,e^{28},&
 de^6&=-2\,e^{17}+e^{24}-e^{35},&\\
 de^7&=2\,e^{16}-e^{25}-e^{34},&
 de^8&=-\sqrt{3}e^{25}+\sqrt{3}\,e^{34}\;.
\end{align*}
Let $G=\LieG{U}(2)$ be the connected subgroup with Lie algebra
\[\lie{g}=\langle e_1,e_6,e_7\rangle\oplus\langle e_8\rangle\cong\su(2)+\lie{u}(1).\]
and fix an invariant connection
\[T=\langle e_2,e_3,e_4,e_5\rangle=\langle t_1,\dotsc,t_4\rangle.\]
If we temporarily identify $T$ with the space of quaternions by
\[1=e_2, i=e_3, j=e_4, k=e_5,\]
then the action of $\lie{g}$ is given by
\[\rho(e_1)=L_i,\; \rho(e_6)=-L_j,\; \rho(e_7)=L_k,\; \rho(\sqrt{3}e_8)=R_k.\]
Notice that $e_8$ was normalized so as to have the same norm as the other elements; consequently, $e_8$ has period $2\pi/\sqrt{3}$.

Let $V=\C^2$. Then we have an equivariant isomorphism $i=\id\colon V\to T$ that induces an element $\beta$ of $\mathcal{L}$, whose components are
\[\beta_1= e_2,\quad\beta_2= e_3,\quad\beta_3= e_4,\quad\beta_4= e_5,\]
and, via the isomorphism $T\cong\Lambda^3T$, another element $\tilde\beta$, given by
\[\tilde\beta_1= e_{345},\quad\tilde\beta_2= -e_{245},\quad\tilde\beta_3= e_{235},\quad\tilde\beta_4= -e_{234}.\]
The isomorphism $i$ also induces an equivariant map
\[j=i\wedge i\colon V\otimes V\to \Lambda^2T,\]
inducing in turn an element $\epsilon$ of $\mathcal{L}$ with components \[\epsilon_k=(a_1\beta_{1}+a_2\beta_{2}+a_3\beta_{3}+a_4\beta_{4})\wedge\beta_k\;.\]
We have two natural invariant contractions $V\otimes V\to\R$, namely the usual scalar product $\cdot$ and a skew-symmetric bilinear form $\sigma$, given by
\[v_1\otimes v_4-v_4\otimes v_1-v_2\otimes v_3+v_3\otimes v_2.\]
The corresponding contractions $m_*(\alpha,\beta)$ will be denoted by $\alpha\beta$, and $\sigma(\alpha,\beta)$ respectively. We shall also indicate the wedge product by juxtaposition, as in $ab\,\sigma(a,b)$ instead of $ab\wedge\sigma(a,b)$.
\begin{table}
 \[\begin{array}{|l|l|}
\hline
p,q & \text{Generators of } \raisebox{0cm}[4mm]{$\Omega^q(V,\Lambda^pT)^{\LieG{U}(2)}$}\\
\hline
0,1 & ab,\sigma(a,b)\\
1,0 & a\beta,\sigma(a,\beta)\\
\hline
0,2 & \sigma(b,b),ab\,\sigma(a,b)\\
1,1 & b\beta,\sigma(b,\beta),ab\,a\beta,ab\,\sigma(a,\beta),\sigma(a,b)\,a\beta,\sigma(a,b)\,\sigma(a,\beta)\\
2,0 & \sigma(\beta,\beta),\sigma(a,\epsilon)\\
\hline
0,3&ab\,\sigma(b,b),\sigma(a,b)\,\sigma(b,b)\\
1,2& ab\,b\beta,ab\,\sigma(b,\beta),\sigma(a,b)\,b\beta,\sigma(a,b)\,\sigma(b,\beta) ,a\beta\,\sigma(b,b),\sigma(a,\beta)\,\sigma(b,b),ab\,\sigma(a,b)\,a\beta ,\\
	 &\hspace{1cm} ab\,\sigma(a,b)\,\sigma(a,\beta)\\
2,1&b\epsilon,\sigma(b,\epsilon),ab\,\sigma(\beta,\beta),ab\,\sigma(a,\epsilon),\sigma(a,b)\,\sigma(\beta,\beta),\sigma(a,b)\,\sigma(a,\epsilon),\sigma(a,\beta)\,b\beta,\sigma(a,\beta)\,\sigma(b,\beta)\\
3,0&a\tilde\beta,\sigma(a,\tilde\beta)\\
\hline
0,4&\sigma(b,b)\,\sigma(b,b)\\
1,3&\sigma(b,b)\,b\beta,\sigma(b,b)\,\sigma(b,\beta),ab\,\sigma(a,b)\,b\beta,ab\,\sigma(a,b)\,\sigma(b,\beta),ab\,a\beta\,\sigma(b,b) ,ab\,\sigma(a,\beta)\,\sigma(b,b)\\
2,2&\sigma(b,b)\,\sigma(\beta,\beta) ,b\beta\,b\beta,b\beta\,\sigma(b,\beta),\sigma(b,\beta)\,\sigma(b,\beta),ab\,b\epsilon, ab\,\sigma(b,\epsilon),\sigma(a,b)\,b\epsilon,\\
	&\hspace{1cm}\sigma(a,b)\,\sigma(b,\epsilon),\sigma(a,\epsilon)\,\sigma(b,b), ab\,\sigma(a,b)\,\sigma(\beta,\beta) ,ab\,\sigma(a,b)\,\sigma(a,\epsilon) ,ab\,\sigma(a,\beta)\,b\beta \\
3,1&b\tilde\beta,\sigma(b,\tilde\beta),ab\,a\tilde\beta ,ab\,\sigma(a,\tilde\beta) ,\sigma(a,b)\,a\tilde\beta ,\sigma(a,b)\,\sigma(a,\tilde\beta)\\
4,0&\beta\tilde\beta\\
\hline
1,4&ab\,\sigma(b,b)\,b\beta ,ab\,\sigma(b,b)\,\sigma(b,\beta) \\
2,3&\sigma(b,b)\,b\epsilon ,\sigma(b,b)\,\sigma(b,\epsilon),ab\,\sigma(a,b)\,b\epsilon ,ab\,\sigma(a,b)\,\sigma(b,\epsilon),ab\,\sigma(b,b)\,\sigma(\beta,\beta) ,ab\,b\beta\,b\beta ,\\
	&\hspace{1cm}ab\,b\beta\,\sigma(b,\beta),\sigma(a,b)\,\sigma(b,b)\,\sigma(\beta,\beta)\\
3,2&b\beta\,b\epsilon,b\beta\,\sigma(b,\epsilon) ,\sigma(b,\beta)\,\sigma(b,\epsilon) ,ab\,b\tilde\beta ,ab\,\sigma(b,\tilde\beta),a\tilde\beta\,\sigma(b,b) ,ab\,\sigma(a,b)\,a\tilde\beta,ab\,\sigma(a,b)\,\sigma(a,\tilde\beta)\\
4,1&ab\,\beta\tilde\beta ,\sigma(a,b)\,\beta\tilde\beta\\
\hline
2,4&\sigma(b,b)\,\sigma(b,b)\,\sigma(\beta,\beta),ab\,\sigma(b,b)\,b\epsilon \\
3,3&\sigma(b,b)\,b\tilde\beta,\sigma(b,b)\,\sigma(b,\tilde\beta),ab\,b\beta\,b\epsilon ,ab\,b\beta\,\sigma(b,\epsilon) ,ab\,\sigma(b,\beta)\,\sigma(b,\epsilon),\sigma(a,b)\,b\beta\,b\epsilon \\
4,2&\sigma(b,b)\,\beta\tilde\beta ,ab\,\sigma(a,b)\,\beta\tilde\beta\\
\hline
3,4&\sigma(b,b)\,b\beta\,b\epsilon ,ab\,\sigma(b,b)\,\sigma(b,\tilde\beta)\\
4,3&b\beta\,b\beta\,b\epsilon ,b\beta\,b\beta\,\sigma(b,\epsilon)\\
\hline
4,4&\raisebox{0cm}[4mm]{$\beta\tilde\beta$}\,\sigma(b,b)\,\sigma(b,b)\\
\hline
\end{array}\]
\caption{The dictionary of  invariant forms on $T\CP^2$.\label{table:GeneratorsTCP2}}
\end{table}
\begin{table}
 \[\begin{array}{|l|l|}
\hline
aa & 2ab\\
\hline
ab & 0\\
\sigma(a,b) & -2 \sigma(a,\epsilon)+\sigma(b,b)-aa\sigma(\beta,\beta)\\
a\beta & b\beta\\
\sigma(a,\beta) & \sigma(b,\beta)\\
\hline
\sigma(b,b) & 2 ab\,\sigma(\beta,\beta)+2 \sigma(a,\beta)\,b\beta+2 \sigma(b,\epsilon)\\
b\beta & 0\\
\sigma(b,\beta) & 0\\
\sigma(\beta,\beta) & 0\\
\sigma(a,\epsilon) & \sigma(a,\beta)\,b\beta+\sigma(b,\epsilon)\\
ab\,\sigma(a,b) & 2 ab\,\sigma(a,\epsilon)-ab\,\sigma(b,b)+aa\,ab\,\sigma(\beta,\beta)\\
ab\,a\beta & - ab\,b\beta\\
ab\,\sigma(a,\beta) & - ab\,\sigma(b,\beta)\\
\sigma(a,b)\,a\beta & -\sigma(a,b)\,b\beta+a\beta\,\sigma(b,b)-2aa\, \sigma(a,\tilde\beta)\\
\sigma(a,b)\,\sigma(a,\beta) & -\sigma(a,b)\,\sigma(b,\beta)+\sigma(a,\beta)\,\sigma(b,b)+2aa\, a\tilde\beta\\
\hline
a\tilde\beta & b\tilde\beta\\
\sigma(a,\tilde\beta) & \sigma(b,\tilde\beta)\\
b\epsilon & - b\beta\,b\beta\\
\sigma(b,\epsilon) & - b\beta\,\sigma(b,\beta)\\
ab\,\sigma(b,b) & -2 ab\,\sigma(a,\beta)\,b\beta-2 ab\,\sigma(b,\epsilon)\\
ab\,b\beta & 0\\
ab\,\sigma(b,\beta) & 0\\
ab\,\sigma(\beta,\beta) & 0\\
ab\,\sigma(a,\epsilon) & -ab\,\sigma(a,\beta)\,b\beta-ab\,\sigma(b,\epsilon)\\
\sigma(a,b)\,\sigma(b,b) &-\frac{3}{2}  aa\, \sigma(b,b)\,\sigma(\beta,\beta)- aa\, b\beta\,b\beta+3 ab\,\sigma(a,b)\,\sigma(\beta,\beta)-2 ab\,b\epsilon\\
&\hspace{2cm}-2 \sigma(a,b)\,\sigma(b,\epsilon)-\sigma(a,\epsilon)\,\sigma(b,b)+\sigma(b,b)\,\sigma(b,b)\\
\sigma(a,b)\,b\beta & -4  aa\, \sigma(b,\tilde\beta)+2 ab\,\sigma(a,\tilde\beta)-2 \sigma(a,b)\,a\tilde\beta+\sigma(b,b)\,b\beta\\
\sigma(a,b)\,\sigma(b,\beta) & 4  aa\, b\tilde\beta-2 ab\,a\tilde\beta-2 \sigma(a,b)\,\sigma(a,\tilde\beta)+\sigma(b,b)\,\sigma(b,\beta)\\
\sigma(a,b)\,\sigma(\beta,\beta) & 3  aa \,\beta\tilde\beta+\sigma(b,b)\,\sigma(\beta,\beta)\\
\sigma(a,b)\,\sigma(a,\epsilon) & -\frac{1}{4}  aa\, \sigma(b,b)\,\sigma(\beta,\beta)-\frac{1}{2}  aa\, b\beta\,b\beta+\frac{1}{2} ab\,\sigma(a,b)\,\sigma(\beta,\beta)-ab\,b\epsilon\\
&\hspace{3cm}-\sigma(a,b)\,\sigma(b,\epsilon)+\frac{3}{2} \sigma(a,\epsilon)\,\sigma(b,b)+\frac{1}{2} aa^2 \beta\tilde\beta \\
a\beta\,\sigma(b,b) & -2  aa\, \sigma(b,\tilde\beta)+6 ab\,\sigma(a,\tilde\beta)-2 \sigma(a,b)\,a\tilde\beta+\sigma(b,b)\,b\beta\\
\sigma(a,\beta)\,\sigma(b,b) & 2  aa\, b\tilde\beta-6 ab\,a\tilde\beta-2 \sigma(a,b)\,\sigma(a,\tilde\beta)+\sigma(b,b)\,\sigma(b,\beta)\\
\sigma(a,\beta)\,b\beta & b\beta\,\sigma(b,\beta)\\
\sigma(a,\beta)\,\sigma(b,\beta) & \sigma(b,\beta)\,\sigma(b,\beta)\\
ab\,\sigma(a,b)\,a\beta & 2  aa\, ab\,\sigma(a,\tilde\beta)+ab\,\sigma(a,b)\,b\beta-ab\,a\beta\,\sigma(b,b)\\
ab\,\sigma(a,b)\,\sigma(a,\beta) & -2  aa\, ab\,a\tilde\beta+ab\,\sigma(a,b)\,\sigma(b,\beta)-ab\,\sigma(a,\beta)\,\sigma(b,b)\\
\hline
\end{array}\]
\caption{Action of $d$ on invariant forms of degree $1$, $2$, $3$.\label{table:TCP2}}
\end{table}

\begin{proposition}
Let $L=\{a,b,\beta,\epsilon\}$, and $B=\{\cdot, \sigma\}$. Then
Algorithm~\ref{algo:} yields a complete set of generators for the vector space of invariant forms.
\end{proposition}
\begin{proof}
Applying Algorithm~\ref{algo:}, we obtain the set of generators appearing in Table~\ref{table:GeneratorsTCP2}. By Theorem~\ref{thm:CompleteDictionary}, it suffices to compute the dimension of $\Lambda^{\LieG{U}(2)}$ and $\Lambda^{\LieG{U}(1)\cdot\LieG{U}(1)}$, where  $\LieG{U}(1)\cdot\LieG{U}(1)$ is the principal stabilizer. The result of this computation is summarized in the following tables:
\[
\begin{array}{|r|r|r|r|}
\hline
 p\backslash q & 0,4    &   1,3&2\\
 \hline
 0,4            & 1   &   0 & 1\\
 1,3            & 0   &   2 &0 \\
 2            & 1   &  0 & 4\\
 \hline
 \end{array}\qquad
\begin{array}{|r|r|r|r|}
\hline
 p\backslash q & 0,4    &   1,3&2\\
 \hline
 0,4            & 1   &   2 & 2\\
 1,3            & 2   &   6 &8 \\
 2            & 2   &  8 & 12\\
 \hline
 \end{array}
\]
Notice that each element appearing in Table~\ref{table:GeneratorsTCP2} vanishes at the origin of $V$ if and only if it contains an $a$ or an $\epsilon$.
Using this fact, it is straightforward to check that Theorem~\ref{thm:CompleteDictionary} is satisfied.
\end{proof}

\subsection{Calabi-Yau cones}
The examples of Section~\ref{sec:TS2} can be regarded either as Calabi-Yau or hyperk\"ahler metrics, as in four dimensions the two notions agree. Not so in eight dimensions: indeed, a Calabi-Yau structure on an $8$-manifold is the special geometry defined by two closed forms $\omega$, $\Psi$ which, at each point, have the form
\begin{equation}
\label{eqn:CY}
\begin{gathered}
\omega=e^{12}+e^{34}+e^{56}+e^{78},\\
\Psi=(e^1+ie^2)\wedge(e^3+ie^4)\wedge(e^5+ie^6)\wedge(e^7+ie^8)\;,
\end{gathered}
\end{equation}
where $e^1,\dotsc,e^8$ is a basis of one-forms. By contrast, a hyperk\"ahler structure on an $8$-manifold is the special geometry defined by three closed forms $\omega_1$, $\omega_2$, $\omega_3$ satisfying
\begin{equation}
\label{eqn:HK}
\omega_1=e^{12}+e^{34}+e^{56}+e^{78},\;\omega_2=e^{13}+e^{42}+e^{57}+e^{86},\;\omega_3=e^{14}+e^{23}+e^{58}+e^{67},
\end{equation}
at each point, where again $e^1,\dotsc,e^8$ is a basis of one-forms. One can relate the two geometries by setting, for example,
\begin{equation}
\label{eqn:CYFromHK}
\omega=\omega_1,\quad \Psi=\frac12(\omega_2\wedge\omega_2-\omega_3\wedge\omega_3)+i\omega_2\wedge\omega_3,
\end{equation}
reflecting the fact that hyperk\"ahler manifolds are Calabi-Yau.

There is a correspondence,  described in \cite{Bar}, between conical metrics with special holonomy and special geometries in one dimension less. In particular, conical hyperk\"ahler structures are related to $3$-Sasaki structures in seven dimensions. Indeed, a $3$-Sasaki structure on a $7$-manifold $M$ can be defined by three contact forms $\eta_1$, $\eta_2$, $\eta_3$ such that the closed forms
\[\omega_i=\eta_i\wedge rdr -\frac12r^2d\eta_i\]
define a conical hyperk\"ahler structure on $M\times\R_+$, where $r$ is a coordinate on $\R_+$.

We now rephrase a well-known result in our language.
\begin{proposition}
\label{prop:3Sasaki}
The one-forms
\begin{equation}
\label{eqn:3Sasaki}
\eta^1=\frac12\sigma(a,b),\quad\eta^2= a\beta, \quad \eta^3=\sigma(a,\beta)
\end{equation}
define a $3$-Sasaki structure on the unit sphere bundle $\mathbb{S}_1\subset T\CP^2$. Accordingly, the complement of the zero section in $T\CP^2$ admits a conical hyperk\"ahler structure, given by
\begin{gather*}
\omega_1=-\frac14\sigma(b,b)+\frac14aa\,\sigma(\beta,\beta)+\frac12\sigma(a,\epsilon),\quad
\omega_2=-\frac12aa^{\frac12}b\beta-\frac12aa^{-\frac12}ab\,a\beta,\\
\omega_3=-\frac12aa^{\frac12}\,\sigma(b,\beta)-\frac12aa^{-\frac12}\,ab\,\sigma(a,\beta),\quad
\end{gather*}
\end{proposition}
\begin{proof}
Let $\mathbb{S}_r$ be the sphere bundle defined by the equation $aa=r^2$. We can extend the $1$-forms $\eta^i$ to the union of the $\mathbb{S}_r$ by
\begin{equation*}
\eta^1=\frac1{2r^2}\sigma(a,b),\quad\eta^2=\frac1r a\beta, \quad \eta^3=    \frac1r\sigma(a,\beta).
\end{equation*}
Notice that these forms are constant along the radial direction; so, we can consistently  identify $\mathbb{S}_1\times\R_+$ with the complement of the zero section in $T\CP^2$. Since the $\eta_i$ are invariant under the homogeneous action of $\SU(3)$, we can work at a point of $\mathbb{S}_r$, where we can assume $a=(r,0,0,0)$. Then
\[\eta_1=\frac1{2r}  b_4,\quad\eta_2=\beta_1,\quad\eta_3=  \beta_4\;.\]
Moreover  $b_1=0$ on $\mathbb{S}_r$. By Table~\ref{table:TCP2}, we obtain
\begin{gather*}
d\eta^1|_{\mathbb{S}_r}=\frac1{r^2}b_{32}+\beta_{23}-2\beta_{14},\quad d\eta^2|_{\mathbb{S}_r}=\frac1r b_3\beta_3+\frac1r b_2\beta_2+\frac1r b_4\beta_4,\\
d\eta^3|_{\mathbb{S}_r}=\frac1r b_3\beta_2+\frac1r\beta_3b_2+\frac1r\beta_1b_4.
\end{gather*}
Thus, the $\eta^i$ are contact forms, and the associated symplectic forms are
\begin{align*}
\omega_1&=-\frac12b_{32}-\frac12r^2\beta_{23}+r^2\beta_{14}-\frac12b_{14}=-\frac14\sigma(b,b)+\frac14aa\,\sigma(\beta,\beta)+\frac12\sigma(a,\epsilon),\\
\omega_2&=-\frac12r(b_3\beta_3+b_2\beta_2+b_4\beta_4)-rb_1\beta_1=-\frac12aa^{\frac12}b\beta-\frac12aa^{-\frac12}ab\,a\beta,\\
\omega_3&=-\frac12r(b_3\beta_2+\beta_3b_2+\beta_1b_4)-rb_1\beta_4=-\frac12aa^{\frac12}\,\sigma(b,\beta)-\frac12aa^{-\frac12}\,ab\,\sigma(a,\beta).
\end{align*}
By construction, or by direct computation using Table~\ref{table:TCP2}, the $\omega_i$ are closed. It only remains to find a local basis of $1$-forms satisfying Equation~\eqref{eqn:HK}; at $a=(r,0,0,0)$, one such basis is given by
\[-\frac{1}{\sqrt2}b_3,\; \frac{1}{\sqrt2}b_2,\; \frac{1}{\sqrt2}r\beta_3,\;\frac{1}{\sqrt2}r\beta_2,\; r\beta_1,\;r\beta_4,\; b_1,\;-\frac1{2}b_4 .\qedhere\]
\end{proof}
The manifold $T\CP^2$ admits other, complete  hyperk\"ahler metrics (see \cite{Calabi:Metriques,BiquardGauduchon:Hyperkahler,DancerSwann:Hyperkahler}), but we shall focus on the Calabi-Yau point of view instead. One can rephrase Proposition~\ref{prop:3Sasaki} by saying that the complement of the zero section in $T\CP^2$ has a three-parameter family of conical symplectic structures
\begin{equation}
\label{eqn:ConicalOmega}
\omega=A\omega_1+B\omega_2+C\omega_3\;,
\end{equation}
and for each  $\omega$ in this family there is a conical Calabi-Yau structure which has $\omega$ as its K\"ahler form. We shall use our list of invariant forms to determine other Calabi-Yau structures (not necessarily conical) defined on a neighbourhood of $\mathbb{S}_1$ in $T\CP^2$ whose K\"ahler form lies in the family \eqref{eqn:ConicalOmega}.

The starting observation is the following: if $M$ is a Calabi-Yau manifold of dimension $8$, each oriented hypersurface inherits an $\SU(3)$-structure, defined by differential forms  $F$, $\Omega$ and $\alpha$ which at each point satisfy
\[F=e^{12}+e^{34}+e^{56},\quad \Omega=(e^1+ie^2)\wedge(e^3+ie^4)\wedge(e^5+ie^6), \quad \alpha=e^7\:;\]
moreover, $F$ and $\alpha\wedge\Omega$ are closed \cite{ContiSalamon}. In the real analytic category, the converse also holds. If one further imposes that the underlying almost-contact structure be contact, i.e. $d\alpha=-2F$, this result can be strengthened as follows:
\begin{theorem}[\cite{ContiFino}]
\label{thm:EmbeddingProperty}
Let $(M,\alpha,F,\Omega)$ be a $7$-manifold with an $\SU(3)$-structure such that
\begin{equation}
 \label{eqn:HypoContact}
d\alpha=-2F,\quad d(\alpha\wedge\Omega)=0.
\end{equation}
Then a neighbourhood of $M\times\{1\}$ in $M\times\R_+$ has a Calabi-Yau structure $(\omega,\Psi)$ such that
\[\omega=\alpha\wedge r dr-\frac12 r^2d\alpha\;.\]
\end{theorem}
In fact, an analogue of this theorem holds in every odd dimension, and five-dimensional examples of structures of this type have been described  in \cite{DeAndres:HypoContact}. Also, solutions of \eqref{eqn:HypoContact} appear to be related to weakly integrable generalized $\Gtwo$-structures \cite{Witt:GeneralizedG2, FinoTomassini}, although there does not seem to be any such structure on $\mathbb{S}_1$.

We can now prove the main result of this section. Notice that we are representing forms on $\mathbb{S}_1$ by pulling back forms on $T^*\CP^2$, so that, say, the radial component of a one-form has no effect on the formulae.
\begin{theorem}
\label{thm:MyBigTheorem}
The seven-dimensional manifold $\mathbb{S}_1$ has a two-parameter family of $\SU(3)$-structures that are solutions of \eqref{eqn:HypoContact} given by
\begin{gather*}
\alpha(B,C) =\frac12\sigma(a,b) + Ba\beta + C\sigma(a,\beta),\\
F(B,C) =  \frac12 \sigma(a,\epsilon)-\frac14\sigma(b,b)+\frac14\sigma(\beta,\beta)-\frac12B b\beta-\frac12C\sigma(b,\beta),\\
\Omega(B,C) = \frac1{4\sqrt{(B^2+C^2)(1+B^2+C^2)}} (\Omega^+ + i\Omega^-),
\end{gather*}
where $B$, $C$ are  constants, at least one of which is not zero, and
\begin{multline*}
\Omega^+= -{2B}\, b\epsilon- {2B(B^{2}+   C^{2}+1 )} \sigma(a,\beta)\,\sigma(b,\beta)
+ {2C} \sigma(a,\beta)\,b\beta\\+ 4  (C^{2}+1)(B^2+C^2)\sigma(a,\tilde\beta)
-B(B^2+C^2) \sigma(a,b)\,\sigma(\beta,\beta)- {2C(1+  B^{2}+ C^{2}) }\sigma(b,\epsilon) \\-{(B^2+C^{2})} \sigma(a,b)\,b\beta
+ 4{CB( B^{2}+ C^{2})} a\tilde\beta\;,
\end{multline*}
\begin{multline*}
\Omega^-= -{(B^2+C^{2})} \sigma(a,b)\,\sigma(b,\beta)+2C \sigma(a,\beta)\,\sigma(b,\beta) -2   B \sigma(b,\epsilon)\\
+2B(1+B^2+C^2)\sigma(a,\beta)\,b\beta- C(B^2+C^2) \sigma(a,b)\,\sigma(\beta,\beta)\\+2C{(1+B^2+C^2)}b\epsilon
- 4CB(B^2+C^2) \sigma(a,\tilde\beta)-4(1+B^2)(B^2+C^2)a\tilde\beta\;.
\end{multline*}
Any two distinct structures in this family have underlying metrics in different conformal classes.
\end{theorem}
\begin{proof}
Equation \eqref{eqn:HypoContact} is proved using Table~\ref{table:TCP2}. To check that the forms actually define an $\SU(3)$-structure, we can work at a point where $a=(1,0,0,0)$. Suppose for the moment that $C\neq 0$. Setting $\delta=B\beta_4-C\beta_1$, $\gamma=2\beta_4-Cb_4$, we compute
\begin{multline*}
\Omega = \frac1{\sqrt{2(B^2+C^2)(1+B^2+C^2)}}(b_2+B\beta_3+C\beta_2 + i(b_{3}+  \beta_{3} C-  \beta_{2} B))\\
\wedge\frac{1}{\sqrt2}(\beta_2- i\beta_3)\wedge
\left(-\frac BC\,\delta+\frac{B^{2}+C^{2}}{2C} \,\gamma + i{(1+B^{2}+C^2)}\delta\right)
\end{multline*}
and
\[F=\frac1C \gamma\wedge\delta +b_{23}-\beta_{23}-B(b_2\beta_2+b_3\beta_3)+C(b_2\beta_3-b_3\beta_2).\]
Thus, the forms $\alpha$, $F$ and $\Omega$ are in standard form with respect to the basis $e^1,\dotsc,e^7$ given by
\begin{gather*}
 e^1=\left(b_2+B\beta_3+C\beta_2\right)/\sqrt2 ,\quad e^2=\left(b_{3}+  \beta_{3} C-  \beta_{2} B\right)/\sqrt2,\\
e^3=\sqrt{(1+B^2+C^2)/2}\,\beta_2,\quad e^4=-\sqrt{(1+B^2+C^2)/2}\,\beta_3,\\ e^5=-\frac{ B}{C\sqrt{(B^2+C^2)(1+B^2+C^2)}}\delta+\frac1{2C}\sqrt{\frac{B^{2}+C^{2}}{(B^2+C^2+1)}}\,\gamma, \\
e^6=\sqrt{\frac{B^2+C^2+1}{B^{2}+C^{2}}}\delta,\quad e^7=\alpha,
\end{gather*}
and we have proved the statement for $C\neq 0$.
Now fix $B\neq 0$, so that $\Omega(B,0)$ can be obtained as the limit
\[\Omega\left(B+\frac 1n,\frac 1n\right)\xrightarrow{n\to+\infty} \Omega(B,0),\]
and consider the associated sequence of coframes constructed as above. It is easy to verify that this sequence converges to a limit coframe $e^1,\dotsc,e^7$; hence, by continuity, the forms $\alpha(B,0)$, $F(B,0)$, $\Omega(B,0)$ define an $\SU(3)$-structure.

Finally, the statement about the underlying Riemannian metric follows from the fact that since $e^1,\dotsc,e^7$ is an orthonormal frame,
\[B=\frac{\langle b_2,\beta_2\rangle}{\langle b_2,b_2\rangle},\quad C=\frac{\langle b_2,\beta_3\rangle}{\langle b_2,b_2\rangle}\;.\qedhere\]
\end{proof}
\begin{remark}
The problem with the case $B=0=C$ is that the map \[(B,C)\to \Omega(B,C)\] is not continuous in the origin. However, we can declare arbitrarily
that
\[\Omega(0,0)=\lim_{n\to +\infty}\Omega\left(\frac1n,0\right)\;,\]
and obtain a valid solution of \eqref{eqn:HypoContact}, given explicitly by
\begin{gather*}
\alpha =\frac12\sigma(a,b),\quad  F =  \frac12 \sigma(a,\epsilon)-\frac14\sigma(b,b)+\frac14\sigma(\beta,\beta)\;,\\
\Omega=- \frac12\sigma(a,\beta)\,\sigma(b,\beta)- \frac12 b\epsilon+\frac12 i\left(- \sigma(b,\epsilon)+ \sigma(a,\beta)\,b\beta\right)\;.
\end{gather*}
This $\SU(3)$-structure is compatible with the $3$\nobreakdash-Einstein-Sasaki structure of Proposition~\ref{prop:3Sasaki}. Indeed, by \eqref{eqn:CYFromHK} a Calabi-Yau structure $(\omega,\Psi)$ is associated to the conical hyperk\"ahler structure, and the forms $F$ and $\alpha\wedge\Omega$ coincide with the restriction to $\mathbb{S}_1$ of the forms $\omega$ and $i\Psi$.
\end{remark}
\begin{remark}
The  $3$\nobreakdash-Einstein-Sasaki structure on $\mathbb{S}_1$ also induces a two-parameter family of solutions of \eqref{eqn:HypoContact}, consisting of $\SU(3)$-structures with the same underlying metric. As a consequence of the last part of Theorem~\ref{thm:MyBigTheorem}, the intersection of the two families consists only of the structure described in the above remark.
\end{remark}
From Theorems~\ref{thm:EmbeddingProperty} and~\ref{thm:MyBigTheorem}, we immediately obtain
\begin{corollary}
For each conical symplectic form $\omega=\omega_1+B\omega_2+C\omega_3$ (with $\omega_i$ as defined in Proposition~\ref{prop:3Sasaki}), there is a neighbourhood of $\mathbb{S}_1$ in $T\CP^2$ carrying a Calabi-Yau metric whose K\"ahler form is $\omega$. The metrics can be chosen  in different conformal classes.
\end{corollary}

\begin{remark}
Considering the complexity of Table~\ref{table:GeneratorsTCP2}, we expect  $T\CP^2$ to  have other special geometries, beside the examples of this section. Indeed, arbitrary choices were involved in the construction of Theorem~\ref{thm:MyBigTheorem}.
\end{remark}

\smallskip
\textbf{Acknowledgements}. This paper is based in part on the author's \emph{Tesi di perfezionamento} \cite{Thesis}, under supervision of S.~Salamon; to him goes the author's gratitude. Thanks are also due to P.~Gauduchon and A.~Swann for their interest and comments on said thesis.
\bibliographystyle{plain}
\bibliography{invariant}

\vskip10pt

\small\noindent Dipartimento di Matematica e Applicazioni, Universit\`a di Milano Bicocca,  Via Cozzi 53, 20125 Milano, Italy.\\
\texttt{diego.conti@unimib.it}
\end{document}